\documentclass[12pt]{scrartcl}

\usepackage{preamble_W}

\newcommand{\setOfReals}{\mathbb{R}}
\newcommand{\setOfNaturals}{\mathbb{N}}
\newcommand{\setOfNonnegativeIntegers}{{\mathbb{N}_0}}
\newcommand{\setOfPositiveReals}{{\setOfReals_{+}}}

\newcommand{\borel}[1]{\mathcal{B} (#1 )}
\newcommand{\spaceOfMeasures}[1]{\mathcal{M}(#1)}
\newcommand{\spaceOfOccupationMeasures}[2]{\mathcal{M}_{#1}(#2)}

\newcommand{\predictableVariation}[1]{ \langle#1\rangle }

\newcommand{\absolute}[1]{| #1 | }

\newcommand{\cadlag}{c\`adl\`ag }

\newcommand{\defeq}{\coloneqq}

\newcommand{\indicator}[2]{\mathsf{1}_{#1}{\left(#2\right)}}
\newcommand{\indic}{\mathsf{1}}

\newcommand{\differential}[1]{\mathrm{d} #1}
\newcommand{\timeDerivative}[1]{\frac{\differential}{\differential t} #1 }

\newcommand{\eqstop}{.}
\newcommand{\eqcomma}{,}

\newcommand{\E}{\mathbb{E}}
\newcommand{\Eof}[1]{\E\left[#1 \right]}

\newcommand{\prob}{\mathbb{P}}
\newcommand{\probOf}[1]{\prob\left(#1\right)}

\newcommand{\history}[1]{\mathcal{F}_{#1}  }

\newcommand{\ConvInProb}{\xrightarrow[]{ \hspace*{4pt}  \text{         P   } }}

\newcommand{\ie}{\textit{i.e.}}
\newcommand{\eg}{\textit{e.g.}}

\newcommand{\nY}{Y^{(n)}}

\newcommand{\nZ}{Z^{(n)}}

\newcommand{\nU}{U^{(n)}}

\newcommand{\nK}{K^{(n)}}

\newcommand{\nM}{M^{(n)}}

\newcommand{\f}{\frac}
\newcommand{\Rt}{\longrightarrow}
\newcommand{\rt}{\rightarrow}
\newcommand{\RT}{\Rightarrow}
\newcommand{\LRT}{\Longrightarrow}

\newcommand{\nrt}{\stackrel{n\rt \infty}\Rt}
\newcommand{\nRT}{\stackrel{n\rt \infty}\LRT}

\newcommand{\prt}{\stackrel{\PP}\Rt}

\newcommand{\leb}{\lambda_{\mbox{Leb}}}

\newcommand{\np}{\noindent}
\newcommand{\non}{\nonumber}
\newcommand{\hs}{\hspace}
\newcommand{\vs}{\vspace}

\newcommand{\om}{\omega}
\newcommand{\Om}{\Omega}

\newcommand{\vep}{\varepsilon}

\newcommand{\s}{\sigma}

\def\SC{\mathcal}

\newcommand{\mfk}{\mathfrak}

\def \triple|{|\! | \! |}
\renewcommand\leq{\ensuremath{\leqslant}}
\def\lf{\left}
\def\ri{\right}
\renewcommand\geq{\ensuremath{\geqslant}}

\def\ot{\otimes}
\def\<{\langle}
\def\>{\rangle}
\def\~{\tilde}
\newcommand{\EE}{\mathbb{E}}
\newcommand{\PP}{\mathbb{P}}

\def\R{\mathbb R}

\newcommand{\gen}{\mathcal{A}}
\newcommand{\fgen}{\mathcal{B}}

\newcommand{\occ}{\Gamma}
\newcommand{\mart}{\SC{M}}
\newcommand{\mar}{\mathsf{M}}
\newcommand{\modu}{\mfk{m}}
\newcommand{\err}{\SC{E}}
\newcommand{\eqb}{z_{C,\star}}

\newcommand{\const}{C}
\numberwithin{equation}{section}

\title{Asymptotic Analysis of the Total Quasi-Steady State Approximation for the Michaelis--Menten Enzyme Kinetic Reactions} 

\author{%
Arnab Ganguly, \emph{Louisiana State University}\\
Wasiur R. KhudaBukhsh\hspace{0.5mm}\href{https://orcid.org/0000-0003-1803-0470}{\includegraphics[width=3mm]{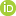}}, \emph{University of Nottingham}\\
}
\date{}

\begin{document}

\maketitle
            
\begin{abstract}
We consider a stochastic  model of the \ac{MM} enzyme kinetic reactions in terms of \acp{SDE}  driven by \acp{PRM}. It has been argued that among various \acp{QSSA} for the deterministic model of such chemical reactions, the \ac{tQSSA} is the most accurate approximation, and it is valid for a wider range of parameter values than the \ac{sQSSA}.  While the \ac{sQSSA} for this model has been rigorously derived from a probabilistic perspective at least as early as 2006 in \cite{Ball:2006:AAM}, a rigorous study of the \ac{tQSSA} for the stochastic model appears missing. We fill in this gap by deriving it as a \ac{FLLN}, and also studying the fluctuations around this approximation as a \ac{FCLT}. 
\end{abstract}



\tableofcontents

\section{Introduction}
    \label{sec:intro}
    The core result of the paper is an intricate stochastic averaging principle that is needed to justify the \acf{tQSSA} of a widely popular model in biochemistry, the \acf{MM} enzyme kinetic reaction system. Despite the surprisingly simple description of the model, a rigorous derivation of the principle via a \acf{FLLN} and a \acf{FCLT}-based analysis of fluctuations providing closed-form expressions of the limits require careful consideration of various functional analytic and probabilistic subtleties. 

\subsection{The \acl{MM} enzyme kinetic reaction system}
    In the simplest form, the \ac{MM} enzyme kinetic reactions  consist of a reversible binding of a substrate and an enzyme into a substrate-enzyme complex, and the production of a product freeing up the bound enzyme \citep{Segel:1975:EK}. Thus, the \ac{MM} reaction network can be schematically represented as follows
    \begin{equation}
    S + E \xrightleftharpoons[\kappa_{-1}]{\kappa_1} C  \xrightharpoonup[]{\kappa_2}    P+ E,
    \label{eq:mm_det}
    \end{equation}
    where $S, E, C$, and $P$ denote molecules of the substrate, the enzyme, the substrate-enzyme complex, and the product. The positive real numbers $\kappa_1$, $\kappa_{-1}$, and $\kappa_2$ are called the reaction rate constants. The time evolution of the concentration of the molecules of $S, E, C$, and $P$ can be described using the following system of \acp{ODE}
    \begin{align}
        \begin{aligned}
        \timeDerivative{S_t} &{} = -\kappa_1 S_t E_t + \kappa_{-1}C_t \eqcomma \quad
        \timeDerivative{E_t} {} = -\kappa_1 S_t E_t + (\kappa_{-1} + \kappa_2) C_t \eqcomma \\
        \timeDerivative{C_t} &{} = \kappa_1 S_t E_t - (\kappa_{-1} + \kappa_2) C_t \eqcomma \quad
        \timeDerivative{P_t} {} = \kappa_2C_t \eqcomma 
        \end{aligned}
        \label{eq:mm_ode}
    \end{align}
    with initial conditions $S_0 = s_0, E_0=e_0, C_0=0, P_0=0$. Here, we have implicitly assumed the law of mass action \citep{Anderson:2011:CTM}. The \ac{MM} system admits two conservation laws: for all $t\geq 0$, 
    \begin{align*} 
        E_0 &{} = e_0 = E_t + C_t  \eqcomma \quad 
        S_0 {} = s_0 = S_t + C_t + P_t. 
    \end{align*}
    Experimental data suggest that the complex $C$ reaches a steady-state rapidly while the species $S, E$, and $P$ remain in their transient states. Therefore, by setting $\timeDerivative{C_t} \approx 0$, we get the steady-state value $C = e_0 S /(K_M+S)$
    where $\kappa_M=(\kappa_{-1}+\kappa_2)/\kappa_1$. The substrate concentration is then given by
    \begin{align}
    \timeDerivative{S_t} &= -\frac{\kappa_2 e_0 S_t}{\kappa_M+S_t}. \label{eq:det_sQSSA_S}
    \end{align}
    This ad hoc approximation is known as the deterministic \ac{sQSSA} of the \ac{MM} enzyme kinetic reaction system in \eqref{eq:mm_det}. The validity of this approximation, controversy due to its rampant misuse, and its many generalizations have been studied in the literature over the last few decades. See \cite{Segel:1988:VSS,Segel:1989:QSS,Stiefenhofer:1998:QSS,Tzafriri:2007:QSS,Tzafriri:2004:TQS,Tzafriri:2003:MMK,Choi2017BeyondMM,Kim2014validity,Kim2020Misuse,Schnell2013Validity,Schnell2000Enzymekinetics,Kang2019QSSA,Eilertsen2024Unreasonable,Pedersen2008tQSSA,Borghans1996QSSA}. It has been argued that the \ac{tQSSA}, which we describe below, is a more accurate approximation for the \ac{MM} system than the \ac{sQSSA} \cite{Kang2019QSSA,Tzafriri:2004:TQS,Kim2020Misuse,Eilertsen2024Unreasonable}. However, even this claim is not without scepticism and controversy. See \cite{Eilertsen2024Unreasonable,Srivastava2025QSSA} for a recent discussion on this topic. Nevertheless, it remains an important mathematical topic that has widespread implications in systems biology, and biochemistry. 

    In sharp contrast to the \ac{sQSSA},  one introduces a new variable $T \defeq  S+C$ in the deterministic \ac{tQSSA}. Note that 
    \begin{align}
    \label{det:tQSSA_eq}
    \begin{aligned}
    \timeDerivative{T} &= -\kappa_2 C,\quad
    \timeDerivative{C} = \kappa_1 \left( \left(T-C\right)\left(e_0-C\right) - \kappa_M C\right),
    \end{aligned}
    \end{align}
    where $\kappa_M=(\kappa_{-1}+\kappa_2)/k_1$ as before. Assume that $\timeDerivative{C} \approx 0$ and $C \le e_0$. Then,  solving the quadratic equation for $C$ and interpreting the smaller root as the `physically meaningful' solution gives the steady-state value of $C$ as 
    \begin{align}
    C &=\frac{\left(e_0+\kappa_M+T\right)-\sqrt{\left(e_0+ \kappa_M+T\right)^2-4e_0 T}}{2}\eqstop \label{tQSSA_C}
    \end{align}
    The time evolution of the new variable $T$ is thus given by 
    \begin{align}
    \timeDerivative{T} &= -\kappa_2 \frac{\left(e_0+\kappa_M+T\right)-\sqrt{\left(e_0+\kappa_M+T\right)^2-4e_0 T}}{2}. \label{tQSSA_T}
    \end{align}
    This heuristic approximation is called the deterministic \ac{tQSSA} of the \ac{MM} enzyme kinetic reaction system in \eqref{eq:mm_det}. We refer the readers to \cite{Eilertsen2024Unreasonable,Kang2019QSSA,Tzafriri:2004:TQS,Pedersen2008tQSSA} for further details on the \ac{tQSSA}.

    \subsection{Our contributions in the context of relevant literature}

While the \ac{sQSSA} has been derived from the stochastic model in \cite{Ball:2006:AAM,Kang2019QSSA,Kang:2014:CLT} and recently in a more general framework in \cite{BharatBaburModel}, a rigorous mathematical derivation of the \ac{tQSSA} from stochastic perspective is notably absent in the literature despite its practical importance. The goal of the paper is to address this gap. Such an endeavor is not merely a matter of mathematical curiosity but is crucial due to its significant practical implications. Specifically, approximations such as the \ac{tQSSA} lead to substantial model reduction facilitating accelerated algorithms for simulations or inference of multiscale \acp{CRN}, for example, see  \cite{HaRa02, RaAr03, JaKr12, GAK15}. A precise mathematical understanding of this approximation including a careful error analysis is essential for assessing the accuracy of these numerical algorithms. For instance, in the important problem of statistical inference, establishing the consistency of estimators of the rate constants obtained from the reduced model relies on a mathematically rigorous derivation of the \ac{tQSSA}.

As mentioned, the highlights of this paper are an \ac{FLLN} (\Cref{thm:tQSSA})  justifying the \ac{tQSSA} for the \ac{MM} enzyme kinetic model and   an \ac{FCLT} (\Cref{thm:Z_V_FCLT}) quantifying the fluctuations around the reduced order limiting model in terms of  an It\^o \ac{SDE} in a suitable scaling regime that encodes the difference in reaction speeds and species abundances. In contrast to a generator-based approach of studying asymptotics of operators and corresponding semigroups in suitable function spaces, we take a more probabilistic route, focusing on a direct analysis of the stochastic process representing the species,  modeled as a system of \acp{SDE} driven by \acp{PRM}. Our framework allows for the limiting averaged model to be  a random \ac{ODE}, rather than the standard deterministic \ac{ODE} as in a typical \ac{tQSSA} approximation.  The limiting behavior of the fast component can roughly be characterized by an invariant distribution associated with the generator of an \ac{ODE}. However,  this \ac{ODE} admits two equilibrium points -- one stable and one unstable -- and hence there are infinitely many such invariant distributions as any convex combination of Dirac measures with atoms at these equilibrium points will be invariant. Although it is intuitively clear that the fast process would converge to the stable equilibrium point, the standard \ac{ODE} theory cannot be applied here to prove this. Instead, establishing that the Dirac measure at the stable equilibrium point is the only invariant distribution governing the eventual steady-state behavior of the fast component requires a careful argument involving asymptotics of a (random) occupation measure encoding both the fast and slow parts of the system, along with certain regularity properties of the equilibrium points.  The \ac{FCLT} (\Cref{thm:Z_V_FCLT}), as is typically the case, requires even more delicate analysis relying on key properties of the solution to a Poisson equation.

Prior rigorous works on stochastic multiscale reaction systems include \cite{Ball:2006:AAM, Crudu2012AAP, Kang:2013:STM, Kang:2014:CLT}. These papers primarily use a generator-based approach building on the results of \cite{Kurtz1992Averaging}, which provides a general framework for stochastic averaging of a class of martingale problems. However, the results from these studies cannot be directly applied to the model considered in this paper. For example, \cite{Crudu2012AAP}, which only focuses on deriving the limiting process but does not include a fluctuation analysis, assumes that the amount of the abundant species and the propensity of fast reactions are of the same order --  a condition that does not hold in our case. Furthermore, it assumes the uniqueness of the invariant distribution for the generator of the fast process, a condition that is common in stochastic averaging related papers and also has been assumed in \cite{Kang:2013:STM, Kang:2014:CLT}, but it does not hold in our case. We also note that while \cite{Kang:2013:STM, Kang:2014:CLT} in theory considers a more general framework, the rigorous verification of the general conditions -- such as the validation of \cite[Condition 2.1 - Condition 2.10]{Kang:2014:CLT} for many \acp{CRN}, including our setup -- is highly non-trivial and would require an extensive, paper-length analysis.

Averaging principles have been extensively studied in  generic settings of continuous diffusion processes (for example, see \cite{ArKoNe06, Neis90, Khas68, Skorokhod1989Asymptotic, Veretennikov1990Averaging, Lip96, PaVe01, PaVe02, PaVe03, FrWe08, FrWe12}), and, in recent years, for systems with jumps; see \cite{HMS15, Xu2017Averaging, BDG18, Shen2022Averaging, Gaee2023Averaging, Mao2024Averaging}. While \acp{CRN} are inherently multiscale and often modeled by jump Markov processes, standard multiscale frameworks with jumps do not directly apply, as they typically assume a clear identification of fast-slow components with fast dynamics modeled by a finite-state Markov chain. In contrast, reaction networks display a more complex hierarchy induced by variability in both reaction speeds and species population levels.  Moreover, such systems are strongly coupled, with both fast and slow reactions, as well as rare and abundant species, influencing individual species' temporal behavior in a highly interconnected manner.

    The rest of the paper is structured as follows. In \Cref{sec:stoch_model}, we describe the stochastic model of the \ac{MM} system, which we describe in terms of \acp{SDE} driven by \acp{PRM}. In \Cref{sec:tQSSA}, we specify the scaling regime for the validity of the \ac{tQSSA}. We provide a rigorous derivation of the \ac{tQSSA} as an \ac{FLLN} in \Cref{sec:FLLN} followed by the \ac{FCLT} in \Cref{sec:FCLT}. 
    Additional mathematical derivations are provided in Appendix~\ref{sec:aux-results}.

    \subsection{Notational conventions} 
    \textbullet\ The space of continuous functions from a metric space $E$ to a metric space $F$ will be denoted by $C(E, F)$ with the subset $C_b(E, F)$ containing the bounded, continuous functions when $F$ is a set of real numbers. 
    The set $D([0,T], F)$ denotes  the space of \cadlag functions from the interval $[0,T]$ to $F$.
    In our case, the spaces $E$ and $F$ will be complete, separable metric spaces. The space $C(E, F)$ will be equipped with the supremum norm, whereas the space $D([0,T], F)$ will be equipped with the Skorokhod topology (see \cite[Chapter 3]{Ethier:1986:MPC}, \cite{Billingsley1999Convergence}, or \cite[Chapter 12]{Whitt2002StochLimits}). \textbullet\ The  Borel $\sigma$-field on a metric space $E$ is denoted by $\borel{E}$. \textbullet\ $\spaceOfMeasures{E}$ will denote the space of finite (non-negative) measures on $E$ equipped with the topology of weak convergence. For $r>0$, $\spaceOfOccupationMeasures{r}{E} \subset \spaceOfMeasures{E}$ will denote the space of (non-negative) measures $\nu$ such that $\nu(E) =r$. \textbullet\ The set of natural numbers, non-negative integers, real numbers, non-negative real numbers are denoted by $\setOfNaturals, \setOfNonnegativeIntegers, \setOfReals$, and $\setOfPositiveReals$ respectively. \textbullet\ The indicator function of a set $A$ will be denoted by $\indicator{A}{\cdot}$, \ie, $\indicator{A}{x} = 1$ if $x\in A$, and zero otherwise. \textbullet\ We use $\leb$ to denote the Lebesgue measure on $\setOfPositiveReals$. \textbullet\ For a \cadlag function $f$, we denote the left-hand limit of the function $f$ at $t$ by $f(t-)$.   For a differentiable function $f:\R^d \rt \R$, the function $\partial_i f$ will denote the first-order derivative of $f$ with respect to the $i$-th coordinate. If the coordinate is clear from the context, we will simply use $\partial f$. The $k$-th order derivative with respect to the $i$-th coordinate (when it exists) will be denoted by $\partial^k_i f$. \textbullet\ $a \wedge b$ and  $a \vee  b$ denote, respectively, the minimum and maximum of two real numbers $a$ and $b$.  \textbullet\ The notation $\stackrel{d}=$ will be used to denote equality in distribution, and $\RT$ will denote weak convergence or convergence in distribution. Convergence in probability with respect to a probability measure $\PP$ will be denoted by $\prt$.
    \textbullet\ Other notations will be introduced when needed.

    \section{Stochastic model}
    \label{sec:stoch_model}
    For each $n \geq 1$, interpreted as a scaling parameter, let  $X_S^{(n)}$, $X_E^{(n)}$, $X_C^{(n)}$, and $X_P^{(n)}$, respectively, denote the species copy numbers of the substrate ($S$), the enzyme ($E$), the enzyme-substrate complex ($C$), and the product ($P$), and  $\kappa_1^{(n)}, \kappa_{-1}^{(n)}$, and $\kappa_2^{(n)}$  the (stochastic) reaction rate constants for the first, the second, and the third reaction.  We model the stochastic process $$X^{(n)} \defeq (X_S^{(n)}, X_E^{(n)}, X_C^{(n)}, X_P^{(n)})$$ as a pure jump Markov process \citep{Ethier:1986:MPC,Norris1997MarkovChains,Bobrowski2005Functional} with 
    the generator $A^{(n)}$ defined by 
    \begin{align*}
        A^{(n)} f(x) &\defeq  \kappa_1^{(n)} x_1 x_2 \left(f(x_1-1, x_2-1, x_3+1, x_4) - f(x) \right) \\
        &{}\quad 
        + \kappa_{-1}^{(n)}x_3 \left( f(x_1+1, x_2+1, x_3-1, x_4) -f(x)  \right) \\ 
        &{}\quad + \kappa_2^{(n)} x_3 \left(f(x_1, x_2, x_3-1, x_4+1) -f(x) \right) \eqcomma 
    \end{align*}
    where $x\defeq (x_1, x_2, x_3, x_4) \in \setOfNonnegativeIntegers^4$ and $f : \setOfNonnegativeIntegers^4 \mapsto \setOfReals$ is bounded measurable function. The functions 
    \begin{align}
        \lambda_1^{(n)}(x) \defeq \kappa_1^{(n)} x_1 x_2\eqcomma \quad \lambda_{-1}^{(n)}(x) \defeq  \kappa_{-1}^{(n)}x_3 \eqcomma \quad \lambda_2^{(n)}(x) \defeq \kappa_2^{(n)}x_3  
    \end{align}
    are often called the propensity or intensity functions associated with the reactions.
    The trajectories of the stochastic process $X^{(n)}$ can be described by means of the following \acp{SDE} (written in the integral form; see \cite{Bremaud2020PointProcess,IkedaWatanabe2014Stochastic})
    \begin{align}
        \begin{aligned}
            X_S^{(n)}(t) &{} = X_S^{(n)}(0) - \int_{[0,\infty) \times [0, t]}\indicator{[0, \lambda_1^{(n)}(X^{(n)} (s-))]}{v} Q_1(\differential{v}\times \differential{s}) \\
            &{} \quad + \int_{[0,\infty) \times [0,t]}\indicator{[0, \lambda_{-1}^{(n)}(X^{(n)} (s-))]}{v} Q_{-1}(\differential{v}\times \differential{s})\eqcomma \\
            X_E^{(n)}(t) &{} = X_E^{(n)}(0) - \int_{[0,\infty)\times [0,t]}\indicator{[0, \lambda_1^{(n)}(X^{(n)} (s-))]}{v} Q_1(\differential{v}\times \differential{s}) \\
            &{} \quad + \int_{[0,\infty)\times [0,t]}\indicator{[0, \lambda_{-1}^{(n)}(X^{(n)} (s-))]}{v} Q_{-1}(\differential{v}\times \differential{s}) \\
            &{} \quad +  \int_{[0,\infty)\times [0,t]}\indicator{[0, \lambda_{2}^{(n)}(X^{(n)} (s-))]}{v} Q_{2}(\differential{v}\times \differential{s})\eqcomma \\
            X_C^{(n)}(t) &{} = X_C^{(n)}(0) + \int_{[0,\infty)\times [0,t]}\indicator{[0, \lambda_1^{(n)}(X^{(n)} (s-))]}{v} Q_1(\differential{v}\times \differential{s}) \\
            &{} \quad - \int_{[0,\infty)\times [0,t]}\indicator{[0, \lambda_{-1}^{(n)}(X^{(n)} (s-))]}{v} Q_{-1}(\differential{v}\times \differential{s}) \\
            &{} \quad -  \int_{[0,\infty)\times [0,t]}\indicator{[0, \lambda_{2}^{(n)}(X^{(n)} (s-))]}{v} Q_{2}(\differential{v}\times \differential{s})\eqcomma \\
            X_P^{(n)}(t) &{} =  X_P^{(n)}(0) + \int_{[0,\infty)\times [0,t]}\indicator{[0, \lambda_{2}^{(n)}(X^{(n)} (s-))]}{v} Q_{2}(\differential{v}\times \differential{s})\eqcomma 
        \end{aligned}
    \end{align}
    where $Q_1, Q_{-1}$, and $Q_2$ are independent \acp{PRM} on $\setOfPositiveReals\times \setOfPositiveReals$ with intensity $\leb\ot\leb$ where $\leb$ is the Lebesgue measure on $\setOfPositiveReals$.  The random  measures $Q_1, Q_{-1}$, and $Q_2$ are defined on the same probability space $(\Omega, \history{}, \prob)$, and are independent of $X^{(n)}(0)$.  We assume $\history{}$ is $\prob$-complete and  associate to $(\Omega, \history{}, \prob)$ the filtration $(\history{t})_{t\ge 0}$ given by 
    \begin{align}
        \history{t} \defeq \sigma\left(X^{(n)}(0), Q_i((0,s]\times A) \mid s \le t, A \in \borel{\setOfPositiveReals}, i=1, -1, 2\right)\eqcomma 
    \end{align}
    for $t>0$ and let $\history{0}$ contain all $\prob$-null sets in $\history{}$. The filtration $(\history{t})_{t\ge 0}$ is right continuous in the sense that 
    \begin{align*}
        \history{t+} \defeq \cap_{s>0}\history{t+s} = \history{t}\eqstop 
    \end{align*}
    Therefore, the filtered probability space $(\Omega, \history{}, (\history{t})_{t\ge 0}, \prob) $ is complete or the \emph{usual  conditions} (see \cite[Definition 2.25]{karatzas1991brownian} or \cite[Definition 1.3]{jacod2003limit}; also called the Dellacherie conditions)  are satisfied. 

    In order to study various averaging phenomena and ensuing \acp{QSSA}, we will consider the scaled stochastic process $$Z^{(n)} \defeq (Z_S^{(n)}, Z_E^{(n)}, Z_C^{(n)}, Z_P^{(n)})$$ where 
    \begin{align}
    \begin{aligned}
        Z_S^{(n)}(t) &{}= n^{-\alpha_S} X_S^{(n)}(n^\gamma t)\eqcomma \quad
        Z_E^{(n)}(t)  = n^{-\alpha_E} X_E^{(n)}(n^\gamma t)\eqcomma  \\
         Z_C^{(n)}(t) &{} = n^{-\alpha_C} X_C^{(n)} (n^\gamma t)\eqcomma \quad
         Z_P^{(n)}(t)  = n^{-\alpha_P} X_P^{(n)}(n^\gamma t), 
    \end{aligned}
        \label{eq:scaled_process}
    \end{align}
     the real numbers $\alpha_S, \alpha_E, \alpha_C$, and $\alpha_P$ are scaling parameters to describe species abundance, and the scaling parameter $\gamma $ is used to speed up or slow down time. In addition to the above scaling parameters, we will also consider scaling parameters $\beta_1, \beta_{-1}$, and $\beta_2$ to describe the speed of the reactions so that we can write 
    \begin{align}
        \kappa_1^{(n)} = n^{\beta_1} \kappa_1\eqcomma \quad \kappa_{-1}^{(n)} = n^{\beta_{-1}} \kappa_{-1}\eqcomma \kappa_2^{(n)} = n^{\beta_2} \kappa_2,
    \end{align}
    for some $n$-free constants $\kappa_1, \kappa_{-1}, \kappa_{2}$. Such parameterizations are common in the stochastic multiscaling literature \citep{Kang:2013:STM,Kang:2014:CLT,Ball:2006:AAM,Kang2019QSSA,BharatBaburModel}. The most common \ac{QSSA} for the \ac{MM} system is the \ac{sQSSA}, which is obtained
    under the following scaling regime:
    \begin{align}
        \begin{aligned}
            & \alpha_S=\alpha_P=1,\quad \alpha_E=\alpha_C=0,\\
           & \beta_1=0,\quad \beta_{-1}=\beta_2=1, 
           &\gamma = 0\eqstop 
        \end{aligned}
            \label{eq:sQSSA_scalings}
    \end{align} 
    as an \ac{FLLN} for the scaled process $\nZ_S$:  
    \begin{align*}
        \nZ_S \ConvInProb Z_S\eqcomma 
    \end{align*}
    as $n\to \infty$, 
    where the limit $Z_S$ lies in $C([0, T], \setOfPositiveReals)$ with probability one and solves \eqref{eq:det_sQSSA_S}. Please see \cite[Theorem 6.1]{Kang:2013:STM} (also \cite{Kang:2014:CLT,Ball:2006:AAM})  or \cite[Theorem 3.1]{BharatBaburModel} for a rigorous derivation of the \ac{sQSSA}.

    \section{The scaling regime for the \acl{tQSSA}}
    \label{sec:tQSSA}
    In order to derive the \ac{tQSSA} directly from the stochastic model described in \Cref{sec:stoch_model}, we assume the same scaling exponents as assumed in \cite{Kang2019QSSA}: 
    \begin{align}
    \begin{aligned}
        & \alpha_S=\alpha_E=\alpha_C=\alpha_P=1, \\
        & \beta_1=\beta_2=0,\quad \beta_{-1}=1,\\
        & \gamma = 0 \eqstop 
    \end{aligned}
    \label{eq:tQSSA_scalings}
    \end{align}
    The interpretation is that all species are abundant, and the binding and the product formation reactions are slower than the unbinding reaction. There is no need to speed up or slow down time. 
    Under the above scaling regime, the trajectories of the scaled stochastic process $Z^{(n)}$ can be described as
    \begin{align}
        \begin{aligned}
            Z_S^{(n)}(t) &{} = Z_S^{(n)}(0) -\frac{1}{n} \int_{[0,\infty)\times [0,t]}\indicator{[0, n^2 \kappa_1 Z_S^{(n)}(s-) Z_E^{(n)}(s-)]}{v} Q_1(\differential{v}\times \differential{s}) \\
            &{} \quad + \frac{1}{n} \int_{[0,\infty)\times [0,t]}\indicator{[0, n^2 \kappa_{-1}Z_C^{(n)} (s-)]}{v} Q_{-1}(\differential{v}\times \differential{s})\eqcomma \\
            Z_E^{(n)}(t) &{} = Z_E^{(n)}(0) -\frac{1}{n} \int_{[0,\infty)\times [0,t]}\indicator{[0, n^2 \kappa_1 Z_S^{(n)}(s-) Z_E^{(n)}(s-)]}{v} Q_1(\differential{v}\times \differential{s}) \\
            &{} \quad + \frac{1}{n} \int_{[0,\infty)\times [0,t]}\indicator{[0, n^2 \kappa_{-1}Z_C^{(n)} (s-)]}{v} Q_{-1}(\differential{v}\times \differential{s}) \\
            &{} \quad + \frac{1}{n} \int_{[0,\infty)\times [0,t]}\indicator{[0, n \kappa_{2}Z_C^{(n)} (s-)]}{v} Q_{2}(\differential{v}\times \differential{s})\eqcomma \\
            Z_C^{(n)}(t) &{} = Z_C^{(n)}(0) +\frac{1}{n} \int_{[0,\infty)\times [0,t]}\indicator{[0, n^2 \kappa_1 Z_S^{(n)}(s-) Z_E^{(n)}(s-)]}{v} Q_1(\differential{v}\times \differential{s}) \\
            &{} \quad - \frac{1}{n} \int_{[0,\infty)\times [0,t]}\indicator{[0, n^2 \kappa_{-1}Z_C^{(n)} (s-)]}{v} Q_{-1}(\differential{v}\times \differential{s}) \\
            &{} \quad - \frac{1}{n} \int_{[0,\infty)\times [0,t]}\indicator{[0, n \kappa_{2}Z_C^{(n)} (s-)]}{v} Q_{2}(\differential{v}\times \differential{s})\eqcomma \\
            Z_P^{(n)}(t) &{} =  Z_P^{(n)}(0) + \frac{1}{n}\int_{[0,\infty)\times [0,t]}\indicator{[0, n \kappa_{2}Z_C^{(n)} (s-)]}{v} Q_{2}(\differential{v}\times \differential{s})\eqstop 
        \end{aligned}
        \label{eq:Z_trajectories}
    \end{align}
    Define the stochastic process $\nZ_V$ by
    \begin{align*}
        \nZ_V(t) \defeq \nZ_S(t)+\nZ_C(t) \text{ for all } t\ge 0 \eqstop 
    \end{align*}
    From \eqref{eq:Z_trajectories}, we immediately see that the process $\nZ_V$ satisfies 
    \begin{align}
        \nZ_V(t) = \nZ_V(0) - \frac{1}{n} \int_{[0,\infty)\times [0,t]} \indicator{[0, n \kappa_2 \nZ_C(s-)]}{v} Q_2(\differential{v}\times  \differential{s})\eqstop 
        \label{eq:Z_V_trajectory}
    \end{align}
    Moreover, observe  that the following two {\em conservation laws} hold:
    \begin{align}
        \begin{aligned}
          \nZ_V(t) + \nZ_P(t) \equiv&\  \nZ_S(t) + \nZ_C(t) + \nZ_P(t)\\
          =&\  \nZ_V(0) + \nZ_P(0) \equiv \ \nK_1 \eqcomma \\
        \nZ_E(t) + \nZ_C(t)  =&\ \nZ_E(0) + \nZ_C(0) \equiv \nK_2 \eqcomma 
        \end{aligned}
        \label{eq:conservation_laws}
    \end{align}
    for some non-negative random variables $\nK_1$, and  $\nK_2$. Notice that the stochastic process $(\nZ_V, \nZ_C)$ itself  is a continuous-time Markov process. Since $\nK_1, \nK_2$ are $\history{0}$-measurable, for any bounded measurable function $f$, the stochastic process 
    $$f(\nZ_V(\cdot), \nZ_C(\cdot))-f(\nZ_V(0), \nZ_C(0))-\int_0^\cdot \gen^{(n)} f(\nZ_v(s),\nZ_C(s)) \differential{s}$$
    is a martingale. Here, the operator $\gen^{(n)}$ is defined by
    \begin{align}
        \begin{aligned}
            \gen^{(n)} f(z_V, z_C) &={} n\kappa_2 z_C \left( f(z_V - \frac{1}{n}, z_C - \frac{1}{n}) - f(z_V, z_C) \right) \\
            &\quad + n^2 \kappa_1 (z_V-z_C)(\nK_2 - z_C)\left(f(z_V, z_C+\frac{1}{n}) - f(z_V, z_C) \right)\\
            &\quad + n^2 \kappa_{-1}z_C \left(f(z_V, z_C- \frac{1}{n}) - f(z_V, z_C) \right)\eqcomma 
        \end{aligned}
        \label{eq:reduced_Z_generator}
    \end{align}
    for measurable functions $f: \setOfPositiveReals^2 \mapsto \setOfReals$.
    It is clear that in the \ac{tQSSA}, the new stochastic process $\nZ_V$ serves as the slow process as opposed to $\nZ_S$ in the \ac{sQSSA}. Clearly, this slow variable depends on the fast component $\nZ_C$, which undergoes rapid jumps at a rate of order $n^2$, creating a strongly coupled system. 
    
    Our goal is to establish a stochastic averaging principle giving a limiting description of $\nZ_V$ as $n\to \infty$. To that extent, define the occupation measure of $(\nZ_C, \nZ_V)$ on the space $(\setOfPositiveReals\times \setOfPositiveReals)\times \setOfPositiveReals$ by 
    \begin{align}
        \occ^{(n)}(A \times B\times [0, t]) \defeq \int_{0}^{t} \indicator{A \times B}{\nZ_C(s), \nZ_V(s)} \differential{s}, 
        \eqcomma 
    \end{align}
    for $A\times B\in \borel{\setOfPositiveReals \times \setOfPositiveReals }$, the set of Borel subsets of $\setOfPositiveReals \times \setOfPositiveReals$, and $t>0$.
    Notice that 
    \begin{align}\label{eq:supp-occ-n}
    \mathrm{supp}(\occ^{(n)}) \subset \{(z_C,z_V) \in \setOfPositiveReals\times \setOfPositiveReals: z_C \leq z_V \}\times \setOfPositiveReals.
    \end{align}
  The occupation measure $\occ^{(n)}$ is a random measure on $(\setOfPositiveReals\times \setOfPositiveReals)\times \setOfPositiveReals$, \ie, an $\spaceOfMeasures{(\setOfPositiveReals\times \setOfPositiveReals)\times\setOfPositiveReals}$-valued random variable. Even though the collection of fast processes $\{\nZ_C : n\ge 1\}$ is not relatively compact because of rapid jumps, the sequence of the occupation measures $\{ \occ^{(n)}: n \ge 1\}$ is well-behaved  and, as we will see in the next section, is  relatively compact. 
  
    
    
    \section{The \acl{FLLN}}
    \label{sec:FLLN}
    In this section, we derive the \ac{tQSSA} as a consequence of the \ac{FLLN} for the scaled process $\nZ_V$ under the scaling regime in \eqref{eq:tQSSA_scalings}. 

    \subsection{Relative compactness}

    \begin{myProposition}
    Assume $\{\nK_1 : n\ge 1\} $ is a tight sequence of random variables. Then,  for any $T>0$, the sequence $\{(\occ^{(n)},\nZ_V) : n\ge 1\}$ is relatively compact as $\spaceOfOccupationMeasures{T}{ \setOfPositiveReals\times  \setOfPositiveReals  \times [0, T]}\times D([0, T], \setOfPositiveReals)$-valued random variables. Furthermore, the limit points of $\nZ_V$ are almost surely in $C([0, T], \setOfPositiveReals)$. 
    \label{prop:rel_compactness} 
    \end{myProposition}
    \begin{proof}[Proof of \Cref{prop:rel_compactness}]
        It is enough to show relative compactness of $\nZ_V$ and $\Gamma_n$ separately. 
        By virtue of \cite[Theorem 2.11]{Budhiraja2019Analysis} (see also \cite[Chapter~4]{Kallenberg2017RandomMeasures}),  the sequence of random measures, $\{\occ^{(n)}: n\ge 1 \}$, is relatively compact if the sequence of the corresponding mean measures (sometimes called intensities) $\{\nu^{(n)}: n\ge 1\}$ given by 
        \begin{align*}
            \nu^{(n)} (A\times B\times [0, t]) & \defeq \Eof{\occ^{(n)}(A\times B\times [0, t]) } \\
            & = \int_{0}^{t} \probOf{\nZ_C(s) \in A, \nZ_V(s) \in B }\differential{s}\eqcomma 
        \end{align*}
        for $A, B \in \borel{\setOfPositiveReals}$,
        is relatively compact. Since $\{\nK_1 : n\ge 1\}$ is relatively compact, for any $\vep>0$, there exists an $R_1\equiv R_1(\vep)>0$ such that
        \begin{align}\label{eq:tight-K-1}
            \inf_n\probOf{0\leq \nK_1 \leq R_1(\vep)}\geq 1-\vep.
        \end{align}
        Now consider the compact set $[0,R_1]\times [0,R_1]\times [0,T]$. Since $\nZ_C(s)\leq \nZ_V(s) \leq \nK_1$ for all $s \ge 0$ by \eqref{eq:conservation_laws}, we have, for all $n\geq 1$,
        \begin{align*}
           \nu^{(n)}([0,R_1]\times [0,R_1]\times [0,T]) \geq T\ \probOf{0\leq \nK_1 \leq R_1} \geq T(1-\vep),
        \end{align*}
        which proves the tightness of the sequence $\{\nu^{(n)}: n\ge 1\}$ and hence, the tightness of $\{\occ^{(n)}: n\ge 1\}$ as a collection of $\spaceOfOccupationMeasures{T}{ \setOfPositiveReals\times  \setOfPositiveReals  \times [0, T]}$-valued random variables.

        We next establish the $C$-tightness of $\nZ_V$ in the space $D([0, T], \setOfPositiveReals)$. See the definition of $C$-tightness in Appendix~\ref{sec:aux-results}. Toward this end, write
        \begin{align}\label{eq:ZV-split}
            \nZ_V(t) & {} = \Phi^{(n)}_V(t)+ \mart^{(n)}_V(t),
        \end{align}
        where
        \begin{align*}
           \Phi^{(n)}_V(t) = &\ \nZ_V(0) - \int_{0}^{t} \kappa_2 \nZ_C(s) \differential{s}  
        \end{align*}
        and the process $\mart^{(n)}_V$ is given by 
        \begin{align*}
            \mart^{(n)}_V(t) \defeq& \int_{0}^{t} \kappa_2 \nZ_C(s) \differential{s} - \frac{1}{n} \int_{[0,\infty)\times [0,t]} \indicator{[0, n \kappa_2 \nZ_C(s-)]}{v} Q_2(\differential{v}\times  \differential{s})\\
            = &\ - \frac{1}{n} \int_{[0,\infty)\times [0,t]} \indicator{[0, n \kappa_2 \nZ_C(s-)]}{v} \tilde Q_2(\differential{v}\times  \differential{s})
        \end{align*}
        is a zero-mean martingale. Now, since $0\leq \sup_{t\leq T}\nZ_V(t) \leq \nK_1$ and  the collection $\{\nK_1 : n \ge 1\}$ is tight, the tightness of the sequence $\{\Phi^{(n)}_V : n\ge 1\}$ in $C([0,T],\R)$ follows from Lemma \ref{lem:int-proc-tight} in Appendix~\ref{sec:aux-results}. 
        Next, we show that $\mart^{(n)}_V \prt 0$ in $D([0,T],\setOfPositiveReals)$ as $n\rt \infty$. In fact, we show the following stronger statement: as $n \rt \infty$, 
        \begin{align}
            \label{eq:mart-conv}
            \sup_{t\leq T}|\mart^{(n)}_V(t)| \prt 0.
        \end{align}
        First observe that $\<\mart^{(n)}_V\>$, the predictable quadratic variation of the martingale $\mart^{(n)}_V$, is given by
        \begin{align*}
           \<\mart^{(n)}_V\>_t =  n^{-1}\int_{0}^{t} \kappa_2 \nZ_C(s) \differential{s} \leq n^{-1}\nK_1t. 
        \end{align*}
        Since $\{\nK_1 : n\ge 1 \}$ is tight by hypothesis, we immediately have 
        \begin{align}
            \label{eq:conv-pred-quadvar}
            \<\mart^{(n)}_V\>_T \prt 0
        \end{align}
        as $n \rt \infty.$
        Next for any positive $\vep, \eta$, the  Lenglart--Rebolledo inequality (see \cite[Lemma 3.7]{Whitt2007MCLT}, and \cite[Remark 4.17]{karatzas1991brownian}) gives us 
        \begin{align}
            \probOf{\sup_{t\leq T}|\mart^{(n)}_V(t)|>\eta}\leq \vep+\probOf{\<\mart^{(n)}_V\>_T > \vep\eta^2}\eqstop 
            \label{eq:Lenglart-Rebolledo}
        \end{align}
        Because of the convergence in \eqref{eq:conv-pred-quadvar}, letting $n\rt \infty$ in the \eqref{eq:Lenglart-Rebolledo} we get 
        \begin{align*}
            \lim_{n\rt \infty}\probOf{\sup_{t\leq T}|\mart^{(n)}_V(t)|>\eta}\leq \vep,
        \end{align*}
        which establishes the convergence in \eqref{eq:mart-conv}. 
       
       This proves that the sequence $\{\nZ_V : n \ge 1\}$ is $C$-tight in the space $D([0,T],\setOfPositiveReals)$ and hence, its limit point(s) lies in  $C([0,T],\setOfPositiveReals).$   
    \end{proof}

    The following lemma about the convergence of certain integrals with respect to the random measure $\occ^{(n)}$ will be crucial in establishing the \ac{FLLN} in \Cref{thm:tQSSA}. 

\begin{myLemma} \label{lem:conv-occ-ZV}
    Assume that $\{\nK_1: n\ge 1\} $ is a tight. 
   \begin{enumerate}
   \item \label{item:int-occ-tight} Let $\phi: \setOfPositiveReals \rt \setOfPositiveReals$ be a non-decreasing continuous function. 
    Then, the sequence of random variables $\lf\{\int_{ \setOfPositiveReals\times  \setOfPositiveReals\times [0, T]}\phi(z_V+z_C)\occ^{(n)}(\differential{z_C}\times \differential{z_V}\times \differential{s})\ri\}$ is tight. Furthermore, as $R \rt \infty$,
    $$\int |\phi(z_V+z_C)|\indic_{\{z_V+z_C>R\}} \occ^{(n)}(\differential{z_C}\times \differential{z_V}\times \differential{s}) \stackrel{\PP} \Rt 0$$
    uniformly in $n$; that is, for any $\eta>0$,
    \begin{align*}
        \sup_n \PP\Bigg( \int & |\phi(z_V+z_C)|\indic_{\{z_V+z_C>R\}} \occ^{(n)}(\differential{z_C}\times \differential{z_V}\times \differential{s}) > \eta\Bigg) \stackrel{R \rt \infty}\Rt 0.
    \end{align*}

    \item \label{item:int-occ-conv} Let $Y, \nY, \ n \geq 1$ be $\R$-valued random variables and 
     $h:\R \times \setOfPositiveReals\times \setOfPositiveReals\times [0,T] \rt \R$,   a continuous function satisfying the following condition: there exists a non-decreasing continuous function $\phi: \setOfPositiveReals \rt \setOfPositiveReals$ such that, for any $y,y' \in \R$, $(z_V,z_C,s) \in \setOfPositiveReals\times  \setOfPositiveReals\times [0,T]$,
   \begin{align*}
   |h(y, z_V,z_C, s)| & \leq \phi(z_V+z_C) |y|, \\
    |h(y, z_V,z_C, s) - h(y', z_V,z_C, s)| &\leq \phi(z_V+z_C) |y-y'|
   \end{align*}
    Suppose that $(\nY, \occ^{(n)}) \RT (Y, \occ)$ as $n \to \infty$. 
    Then, for $t>0$,  as $n\rt \infty$, 
    \begin{align}\label{eq:occ-int}
    \begin{aligned}
        &\int_0^t h(\nY, \nZ_V(s),\nZ_C(s),s)\differential{s}\\
        &\equiv\  \int_{ \setOfPositiveReals\times  \setOfPositiveReals\times [0, t]}h(\nY, z_V,z_C,s) \occ^{(n)}(\differential{z_C}\times \differential{z_V}\times \differential{s})\\
       &\  \RT \int_{\setOfPositiveReals\times  \setOfPositiveReals \times [0, t] }h(Y, z_V,z_C,s) \occ(\differential{z_C}\times \differential{z_V}\times \differential{s}).
       \end{aligned}
    \end{align}
 If $(\nY, \occ^{(n)}) \rt (Y,\occ)$ a.s. or in probability,  then the convergence in  \eqref{eq:occ-int} holds in probability. 
 \end{enumerate}
\end{myLemma}
\begin{proof}[Proof of \Cref{lem:conv-occ-ZV}]

\eqref{item:int-occ-tight} Note that the conservation laws in \eqref{eq:conservation_laws} imply that $\nZ_V(t)+\nZ_C(t) \leq 2\nK_1$; hence, by the assumption that $\phi$ is non-decreasing
\begin{align*}
    \int_{ \setOfPositiveReals\times  \setOfPositiveReals \times [0, T]}\phi(z_V+z_C)\occ^{(n)}(\differential{z_C}\times \differential{z_V}\times \differential{s})
    & =\ \int_{[0, T]}\phi(Z^{(n)}_V(s)+Z^{(n)}_C(s)) \differential{s}\\
    &\leq\ \phi(2\nK_1) T.
\end{align*}
 The assertion now follows from the tightness of $\{\phi(2\nK_1): n\ge 1\}$, which holds since $\{\nK_1 : n \ge 1\}$ is tight by hypothesis and the function $\phi$ is continuous (see \cite[Lemma 3.1]{Whitt2007MCLT}).

For the next part, fix  $\eta>0$ and $\vep>0$. Now, by the tightness of $\{\nK_1 : n \ge 1\}$, choose an $R_\vep$ such that
    $\sup_n \probOf{\nK_1>R_\vep/2} \leq \vep.$
Then, for any $R \geq R_\vep$, the fact that $\phi$ is non-decreasing gives the inequality
\begin{align*}
    \int \phi(z_V+z_C)\indic_{\{z_V+z_C>R\}} \occ^{(n)}(\differential{z_C}\times \differential{z_V}\times \differential{s})
    \leq \phi(2\nK_1) \indic_{\{\nK_1 > R/2\}}, 
\end{align*}
which shows that
\begin{align*}
    \sup_n \PP\Bigg( \int & \phi(z_V+z_C)\indic_{\{z_V+z_C>R\}} \occ^{(n)}(\differential{z_C}\times \differential{z_V}\times \differential{s}) > \eta\Bigg) \\
    &\leq\ \sup_n \PP\lf(\phi(2\nK_1) \indic_{\{\nK_1 > R/2\}}> \eta\ri) \leq\ \sup_n\PP(\nK_1 > R/2) \leq\ \vep.
\end{align*}

\np
\eqref{item:int-occ-conv} 
By the Skorohod representation theorem \citep[Theorem 5.31]{Kallenberg2021Foundations}, there exist a probability space $(\tilde \Omega, \tilde{\SC{F}}, \tilde{\PP})$ and random variables  $(\tilde{Y}, \tilde{\occ})$, $(\tilde{Y}^{(n)}, \tilde{\occ}^{(n)}), n \geq 1$ defined on this space such that
   \begin{align*}
      (\tilde{Y}, \tilde{\occ}) &\stackrel{d} = (Y, \occ), \quad  
      (\tilde{Y}^{(n)}, \tilde{\occ}^{(n)}) \stackrel{d} = (\nY, \occ^{(n)}), \quad 
      (\tilde{Y}^{(n)}, \tilde{\occ}^{(n)}) \nrt (\tilde{Y}, \tilde{\occ}), \quad \text{a.s.}
   \end{align*} 
Since $ \int h(\nY, \cdot, \cdot, \cdot) \differential{\occ^{(n)}}
   \stackrel{d} =  \int h(\tilde{Y}^{(n)}, \cdot, \cdot, \cdot) \differential{\tilde{\occ}^{(n)}}$ and similarly,  $ \int h(Y, \cdot, \cdot, \cdot) \differential{\occ}
   \stackrel{d} =  \int h(\tilde{Y}, \cdot, \cdot, \cdot) \differential{\tilde{\occ}}$, 
the assertion follows once we show 
that the convergence in \eqref{eq:occ-int} holds in probability with $\nY, Y, \occ^{(n)}, \occ$ replaced by $\tilde{Y}^{(n)}, \tilde{Y}, \tilde{\occ}^{(n)}, \tilde{\occ}.$ 
To this end, write the integrand as
$h(\tilde{Y}^{(n)}, \cdot, \cdot, \cdot) = h(\tilde{Y}, \cdot, \cdot, \cdot) + (h(\tilde{Y}^{(n)}, \cdot, \cdot, \cdot) - h(\tilde{Y}, \cdot, \cdot, \cdot))$
and notice that 
\begin{align*}
\Bigg|\int_{ \setOfPositiveReals\times  \setOfPositiveReals \times [0, t]} & (h(\tilde{Y}^{(n)}, z_V,z_C,s) - h(\tilde{Y}, z_V,z_C,s))\tilde{\occ}^{(n)}(\differential{z_C}\times \differential{z_V}\times \differential{s})\Bigg|\\
 &\leq \  |\tilde{Y}^{(n)} - \tilde{Y}|\int_{  \setOfPositiveReals\times  \setOfPositiveReals \times [0, t] } \phi(z_V+z_C)\tilde{\occ}^{(n)}(\differential{z_C}\times \differential{z_V}\times \differential{s})\\
 &\stackrel{\tilde \PP} \Rt \ 0, \quad \text{as } {n\rt \infty}
\end{align*}
since $\lf\{\int \phi(z_V+z_C)\tilde{\occ}^{(n)}(\differential{z_C}\times \differential{z_V}\times \differential{s}) \stackrel{d}= \int\phi(z_V+z_C)\occ^{(n)}(\differential{z_C}\times \differential{z_V}\times \differential{s})\ri\}$ is tight by \eqref{item:int-occ-tight}. 
Thus, the assertion will be proved once we show that  as $n \rt \infty$,
   \begin{align} \label{eq:prob-conv-occ-h-0}
    \int_{ \setOfPositiveReals\times  \setOfPositiveReals \times [0, t]}  h(\tilde Y, z_V,z_C,s)  (\tilde{\occ}^{(n)} -\tilde{\occ})(\differential{z_C}\times \differential{z_V}\times \differential{s}) \stackrel{\tilde \PP} \Rt 0,
   \end{align} 
    that is, for any fixed $\eta>0$, $\vep>0$, there exist an $n_0>0$ such that for all $n>n_0$, we have
    \begin{align}
        \label{eq:prob-conv-occ-h}
        \tilde\PP\lf(\lf|\int_{ \setOfPositiveReals\times  \setOfPositiveReals \times [0, t]}  h(\tilde Y, z_V,z_C,s)  (\tilde{\occ}^{(n)} -\tilde{\occ})(\differential{z_C}\times \differential{z_V}\times \differential{s})\ri| > \eta\ri) \leq \vep.
    \end{align}
To this end, first observe that  the assumption on $h$ implies 
   \begin{align} \label{eq:h-int-fin}
   \int_{ \setOfPositiveReals\times  \setOfPositiveReals \times [0, T]}|h(\tilde Y, z_V,z_C,s)| \tilde{\occ}(\differential{z_C}\times \differential{z_V}\times \differential{s}) < \infty, \quad \mbox{a.s.}
   \end{align}
   This requires some argument. Since $\tilde{\occ}^{(n)} \nrt \tilde{\occ}$ a.s. in $\spaceOfOccupationMeasures{T}{ \setOfPositiveReals\times  \setOfPositiveReals \times [0, T]}$ (topologized by weak convergence), and $\phi$ is continuous,  by a version of Fatou's lemma \cite[Equation 1.5]{FeKaZa13},
   \begin{align*}
 \int_{ \setOfPositiveReals\times  \setOfPositiveReals \times [0, T]} & |h(\tilde Y,z_V,z_C,s)| \tilde\occ(\differential{z_C}\times \differential{z_V}\times \differential{s})\\
 \leq&\ |\tilde Y|  \int_{\setOfPositiveReals\times  \setOfPositiveReals \times [0, T]} \phi(z_V+z_C)\tilde\occ(\differential{z_C}\times \differential{z_V}\times \differential{s})\\
 \leq&\ |\tilde Y| \liminf_{n\rt \infty}\int_{[0, T]\times \setOfPositiveReals\times  \setOfPositiveReals}\phi(z_V+z_C)\tilde{\occ}^{(n)}(\differential{z_C}\times \differential{z_V}\times \differential{s})\\
 < &\ \infty, \quad \mbox{a.s.}
   \end{align*}
The last inequality is a consequence of the tightness of the sequence of random variables  $\lf\{\int \phi(z_V+z_C)\tilde{\occ}^{(n)}(\differential{z_C}\times \differential{z_V}\times \differential{s}): n \ge 1 \ri\}$ and the standard fact that if a sequence of random variables $\{V^{(n)}: n \ge 1\}$ is tight, then $\liminf_{n\rt \infty}|V^{(n)}| < \infty,$ a.s.

For $R>0$, define the bounded set 
    $$\SC{B}_{+,R} \equiv \{(z_V,z_C) \in \setOfPositiveReals\times  \setOfPositiveReals: z_V+z_C  < R \},$$
    and denote its closure by $\bar{\SC{B}}_{+,R}.$ Now, \eqref{eq:h-int-fin} implies that 
\begin{align}\label{eq:h-int-tail}
\int_{\bar{\SC{B}}^c_{+,R} \times [0, T]} | h(\tilde Y, z_V,z_C,s)| \tilde\occ(\differential{z_C}\times \differential{z_V}\times \differential{s}) \stackrel{R \rt \infty}\Rt 0 \quad \mbox{a.s.}
\end{align}

 Next, for any $R>0$,   by Urysohn's Lemma, \cite[Page 122]{foll99},  there exists a continuous function $\eta_{R_1}: \setOfPositiveReals\times  \setOfPositiveReals \rt [0, 1]$ such that
    $\eta_{R_1} \equiv 1$ on $\bar{\SC{B}}_{+,R}$ and $\eta_{R} \equiv 0$ on $\SC{B}^c_{+,R+1}$. For a fixed $y\in \R$, define $\tilde h_{R}(y,\cdot) \in C(\setOfPositiveReals\times  \setOfPositiveReals  \times [0, T], \R)$ by 
    $$\tilde h_{R}(y,z_V,z_C,t)  \defeq h(y,z_V,z_C,t)\eta_{R}(z_V,z_C).$$
   Notice that by construction, for any fixed $y \in \R$, $\tilde h_{R}(y,\cdot)$ is actually compactly supported (hence bounded)
 and $|\tilde h_{R}(y,z_V,z_C,t)| \leq |h(y,z_V,z_C,t)|$.
    Write
    \begin{align}\label{eq:int-diff-h-split}
       \begin{aligned}
     &\int_{ \setOfPositiveReals\times  \setOfPositiveReals \times [0, T] }  h(\tilde Y,z_V,z_C,s)  (\tilde{\occ}^{(n)} -\tilde{\occ})(\differential{z_C}\times \differential{z_V}\times \differential{s})\\
     & =\ \int_{\setOfPositiveReals\times  \setOfPositiveReals \times [0, T]}\tilde h_{R}(\tilde Y,z_V,z_C,s) (\tilde{\occ}^{(n)} -\tilde{\occ})(\differential{z_C}\times \differential{z_V}\times \differential{s})  +  \tilde\err^{(n)}_{h,R}(T),
       \end{aligned}
    \end{align}
    where the error term $\tilde\err^{(n)}_{h,R}(T)$ can be estimated as follows:
    \begin{align}\label{eq:err-ineq}
    \begin{aligned}
       |\tilde\err^{(n)}_{h,R}(T)| \leq&\ 2\int_{ \bar{\SC{B}}^c_{+,R} \times [0, T]}|h(\tilde Y, z_V,z_C, s)| (\tilde{\occ}^{(n)} + \tilde{\occ})(\differential{z_C}\times \differential{z_V}\times \differential{s})\\
       \leq &\ 2|\tilde Y| \int \phi(z_V+z_C)\indic_{\{z_V+z_C>R\}} \tilde{\occ}^{(n)}(\differential{z_C}\times \differential{z_V}\times \differential{s}) \\
       &{}\quad 
       +  2\int_{\bar{\SC{B}}^c_{+,R} \times [0, T]}|h(\tilde Y, z_V,z_C, s)|  \tilde{\occ}(\differential{z_C}\times \differential{z_V}\times \differential{s}).
    \end{aligned}
    \end{align}

 We  claim that as $R \rt \infty$, $|\tilde\err^{(n)}_{h,R}(T)| \stackrel{\tilde \PP} \Rt 0$ uniformly in $n$; that is, for any $\eta>0, \vep>0$,
\begin{align}\label{eq:err-unif-conv-main}
   \sup_n \tilde\PP\lf(|\tilde\err^{(n)}_{h,R}(T)|>\eta\ri)\stackrel{R \rt \infty}\Rt 0.
\end{align}
Since the second term in \eqref{eq:err-ineq} converges to $0$ a.s. as $R \rt \infty$ by \eqref{eq:h-int-tail}, we only need to show that as $R \rt \infty$
\begin{align}\label{eq:err-unif-conv}
|\tilde Y| \int \phi(z_V+z_C)\indic_{\{z_V+z_C>R \}} \tilde{\occ}^{(n)}(\differential{z_C}\times \differential{z_V}\times \differential{s}) \stackrel{\tilde \PP} \Rt 0, 
\end{align}
uniformly in $n$. 

Since $\int \phi(z_V+z_C)\indic_{\{z_V+z_C>R\}} \differential{\tilde{\occ}^{(n)}} \stackrel{d}= \int \phi(z_V+z_C)\indic_{\{z_V+z_C>R\}} \differential{\occ^{(n)}}$, 
it follows by \eqref{item:int-occ-tight} that  as $R \rt \infty$, it converges to $0$ uniformly in $n$. Now let $J_\vep$ be such that $\tilde\PP (|\tilde Y| >J_\vep) \leq \vep/2$, and then choose $R_\vep$ such that   
\begin{align*}
        \sup_n \tilde\PP\Bigg( \int & |\phi(z_V+z_C)|\indic_{\{z_V+z_C>R\}} \tilde{\occ}^{(n)}(\differential{z_C}\times \differential{z_V}\times \differential{s}) > \eta/J_\vep\Bigg) \leq \vep/2.
\end{align*}
It follows that 
\begin{align*}
        \sup_n  \tilde\PP\Bigg( & |\tilde Y|\int |\phi(z_V+z_C)|\indic_{\{z_V+z_C>R\}} \tilde{\occ}^{(n)}(\differential{z_C}\times \differential{z_V}\times \differential{s}) > \eta\Bigg)\\
        \leq &\ \sup_n \tilde\PP\Bigg( \int |\phi(z_V+z_C)|\indic_{\{z_V+z_C>R\}} \tilde{\occ}^{(n)}(\differential{z_C}\times \differential{z_V}\times \differential{s}) > \eta/J_\vep\Bigg) \\
        & \quad + \tilde\PP (|\tilde Y| >J_\vep)\ \leq \  \vep.
\end{align*}
This proves the convergence in \eqref{eq:err-unif-conv} and hence the claim.

 Next,  for almost all $\tilde \om \in \tilde \Om$, $\occ^{(n)}(\tilde \om) \rt \occ(\tilde \om)$ (in the weak-convergence topology of $\spaceOfOccupationMeasures{T}{\setOfPositiveReals\times  \setOfPositiveReals \times [0, T]}$), and since the mapping $(z_V,z_C, s) \rt \tilde h_{R}(\tilde Y(\tilde \om),z_V,z_C,t)$ is continuous and bounded for any $R>0$, we have 
\begin{align}
    \label{eq:weak-conv-occ}
    \int  \tilde h_{R}(\tilde Y, z_V,z_C,s)  (\tilde{\occ}^{(n)} -\tilde\occ)(\differential{z_C}\times \differential{z_V}\times \differential{s}) \nrt 0 \text{ a.s.}
\end{align}
Because of this and \eqref{eq:err-unif-conv-main}, \eqref{eq:int-diff-h-split} readily implies \eqref{eq:prob-conv-occ-h-0}. Indeed, 
first choose $R_1 \equiv R_1(\vep)$ such that 
\begin{align*}
   \sup_n \tilde\PP\lf(|\tilde\err^{(n)}_{h,R_1}(T)|>\eta/2\ri) \leq \vep/2.
\end{align*}
Since almost sure convergence in \eqref{eq:weak-conv-occ} implies convergence in probability,  choose $n_0$ such that for all $n>n_0$, we have 
         \begin{align}
        \label{eq:prob-conv-occ-h-trunc}
        \tilde\PP\lf(\lf|\int_{ \setOfPositiveReals\times  \setOfPositiveReals \times [0, t] }  \tilde h_{R_1}(\tilde Y, z_V,z_C,s)  (\tilde{\occ}^{(n)} -\tilde{\occ})(\differential{z_C}\times \differential{z_V}\times \differential{s})\ri| > \eta/2\ri) \leq \vep/2.
    \end{align}   
 The claim in  \eqref{eq:prob-conv-occ-h} now follows from \eqref{eq:int-diff-h-split} and \eqref{eq:prob-conv-occ-h-trunc}.
\end{proof}

Before we state and prove the \ac{FLLN} in \Cref{thm:tQSSA}, we need to introduce a first order differential operator, which will be crucial in the proof.
 
\subsection{A first order differential operator}
For a constant $K$ and a fixed $z_V \in \setOfPositiveReals$, define the operator $\fgen_{K,z_V}$ by
        \begin{align}
            \fgen_{K,z_V} g(z_C) \defeq \left(\kappa_1 (z_V-z_C)(K - z_C) - \kappa_{-1}z_C \right) \partial g(z_C)\eqcomma 
            \label{eq:fast_gen}
        \end{align}
    for $g \in C^{(2)}_{c}(\setOfPositiveReals)$, the space of twice continuously differentiable functions with compact support in $\setOfPositiveReals$. The operator $\fgen_{K_2,z_V}$ can be interpreted as the limiting generator of the fast process when the slow component is frozen at  state $z_V$ and $K_2$ is a limit point of $\nK_2.$
    
    Clearly, for a fixed $K$ and $z_V$, the operator $\fgen_{K,z_V}$ is the infinitesimal generator of the \ac{ODE}:
\begin{align}
\label{eq:ode-eqb-eqn}
    \frac{\differential{z_C(t)}}{\differential{t}}= \kappa_1 (z_V-z_C(t))(K - z_C(t)) - \kappa_{-1}z_C(t).
\end{align}
This \ac{ODE} has two equilibrium points, $\eqb^{-}(K,z_V), \eqb^{+}(K,z_V)$, which are the roots of the quadratic equation (in variable $z_C$)
$$\kappa_1 (z_V-z_C)(K - z_C) - \kappa_{-1}z_C=0\eqcomma $$
given by
\begin{align}
    \label{eq:ode-eqb-pts}
    \eqb^{\pm}(K,z_V)=&\  \frac{(z_V + K)\kappa_1 + \kappa_{-1}  \pm \sqrt{ \left((z_V + K)\kappa_1 + \kappa_{-1}\right)^2 - 4 \kappa_1^2 z_V K } }{ 2\kappa_1}.
\end{align}
Therefore,  the number of invariant distributions corresponding to the operator, $ \fgen_{K,z_V}$, is infinite. Indeed,  for any $0\leq \alpha\leq 1$, the probability measure $\alpha\delta_{\eqb^{-}(K,z_V)} +(1-\alpha)\delta_{\eqb^{+}(K,z_V)}$ is an invariant distribution of $ \fgen_{K,z_V}$.

We now make a few observations about the equilibrium points that are crucial for our main results. First,  it is easy to see that $\eqb^{-}(K,z_V)$ is a stable equilibrium, while $\eqb^{+}(K,z_V)$ is unstable. 
 Next note that we have an immediate lower estimate of the discriminant,
\begin{align}\label{eq:disc}
\begin{aligned}
D(K,z_V) \equiv&\ \left((z_V + K)\kappa_1 + \kappa_{-1}\right)^2 - 4 \kappa_1^2z_V K\\
=&\ (z_V-K)^2\kappa_1^2+2\kappa_{-1}(z_V+K)\kappa_1+\kappa^2_{-1}\\
\geq&\  (z_V-K)^2\kappa_1^2 \vee \kappa_{-1}^2\eqstop 
\end{aligned}
\end{align}
Consequently, we have 
\begin{align*}
    \eqb^+(K,z_V) - z_V =&\ \f{z_V+K}{2}+\f{\kappa_{-1}}{2\kappa_1}+\f{1}{2\kappa_1}\sqrt{D(K,z_V)} - z_V\\
    \geq &\ \f{z_V+K}{2}+\f{\kappa_{-1}}{2\kappa_1} +\f{z_V-K}{2}-z_V\ \geq\ \f{\kappa_{-1}}{2\kappa_1}.
\end{align*}
Thus, the following lower bound on the distance between the unstable equilibrium point $\eqb^{+}(K,z_V)$ and $z_V$ is immediately obtained:
\begin{align}
    \label{eq:min-dist-eqb}
    \min_{z_C \leq z_V}\{\eqb^+(K,z_V) - z_C\} = \eqb^+(z_V) - z_V \geq \f{\kappa_{-1}}{2\kappa_1}
\end{align}
In particular, \eqref{eq:min-dist-eqb} shows that the unstable equilibrium point $\eqb^{+}(K,z_V) \in (z_V, \infty)$ for $\kappa_1, \kappa_{-1} \in (0,\infty).$ It is also easy to see by simple algebra that the stable equilibrium point $\eqb^{-}(K,z_V) \in (0,z_V]$ for any  $\kappa_1>0, \kappa_{-1} \geq 0$. Indeed, using the observation that $z_V+K\geq |K - z_V|$ we see from \eqref{eq:disc} that $D(K,z_V) \geq |K-z_V|\kappa_1+\kappa_{-1}$,  which gives $\eqb^{-}(K,z_V) \leq z_V \wedge K.$ The notation $a \wedge b$ denotes the minimum of $a$ and $b$.

Next, notice that the mapping $z_V \mapsto \eqb^{\pm}(K, z_V)$ is differentiable (with respect to $z_V$), and the derivative is given by
\begin{align}
    \label{eq:deriv-eqb}
    \partial\eqb^{\pm}(K,z_V)=&\ \f{1}{2} \pm \f{(z_V + K)\kappa_1 + \kappa_{-1}-2\kappa_1K}{\sqrt{ \left((z_V + K)\kappa_1 + \kappa_{-1}\right)^2 - 4 \kappa_1^2 z_V K }}.
\end{align}
Therefore, by virtue of  \eqref{eq:disc}, the mapping $z_V \Rt \partial\eqb^{\pm}(K,z_V) $ has at most linear growth, that is, for some constant $\const_{*,1}>0$, 
\begin{align}
    \label{eq:eqb-deriv-lin}
    | \partial\eqb^{\pm}(K,z_V)| \leq \const_{*,1}(1+z_V), \quad z_V \in \setOfPositiveReals.
\end{align}
Similar algebra also shows that for some constant  $\const_{*,2}>0$, 
\begin{align}
    \label{eq:eqb-deriv-lin-2}
    | \partial^2\eqb^{\pm}(K,z_V)| \leq \const_{*,2}(1+z^2_V), \quad z_V \in \setOfPositiveReals.
\end{align}

The following easy lemma, whose proof is included for the sake of completeness, identifies the invariant distribution of the operator $\fgen_{K,z_V}$ that is supported on $[0,z_V]$.
\begin{myLemma}
    \label{lem:fast-inv-dist}
 Let $K>0, z_V>0$.    Let $\mu$ be a probability measure supported on $[0,z_V]$ satisfying
 \begin{align*}
     \int_{\setOfPositiveReals}\fgen_{K,z_V}g(z_C)\mu(\differential{}z_C)=0\eqcomma 
 \end{align*}
 for any  $g \in C^{(2)}_c(\setOfPositiveReals, \R)$, the space of twice continuously differentiable functions with compact support.
 Then $\mu = \delta_{\eqb^{-}(K,z_V)}$.
\end{myLemma}
\begin{proof}
Choose $g \in C^{(2)}_c(\setOfPositiveReals, \R)$ such that
$$g(z_C) = \int_0^{z_C}\left(\kappa_1 (z_V-u)(K - u) - \kappa_{-1}u \right)\differential{u}, \qquad \mbox{for } z_C \leq z_V.$$  
Since $ \mathrm{supp}(\mu) \subset [0,z_V]$, we have 
\begin{align*}
 \int_{\setOfPositiveReals}\fgen_{K,z_V}g(z_C)  \mu(\differential{}z_C)
 =&\  \int_{\R_{+}} \left(\kappa_1 (z_V-z_C)(K - z_C) - \kappa_{-1}z_C \right)\partial g(z_C) \mu(\differential{}z_C)\\
  =&\ \int_{0}^{z_V} \left(\kappa_1 (z_V-z_C)(K - z_C) - \kappa_{-1}z_C \right)^2  \mu(\differential{}z_C)\\
  =&\ 0, 
\end{align*}
which shows that for $\mu$-a.a $z_C \in [0,z_V]$ 
\begin{align*}
    \kappa_1 (z_V-z_C)(K - z_C) - \kappa_{-1}z_C =0.
\end{align*}
But from the discussion above the only $z_C$ in $[0,z_V]$ satisfying the above equation is given by $z_C = \eqb^{-}(K_2,z_V)$ (see \eqref{eq:ode-eqb-pts}). Therefore, we must have $\mu=\delta_{\eqb^{-}(K,z_V)}$.
\end{proof}

We are now ready to state and prove the main result in this section, which establishes the \ac{tQSSA} for the \ac{MM} reaction system as an \ac{FLLN}.

    \begin{myTheorem}
        Assume that the sequence of random variables $\{\nK_1\}$ is tight, and as $n\rt \infty$,  $(\nK_2, \nZ_V(0)) \RT (K_2, Z_V(0))$ where the limit $(K_2, Z_V(0))$ can be random. Then, as $n \rt \infty$, $\nZ_V \RT Z_V$, where the path-space of $Z_V$ is $C([0,T],\setOfPositiveReals)$, and for $\PP$-a.a. $\om \in \Omega$, $Z_V \equiv Z_V(\cdot, \om)$ solves the random \ac{ODE}
        \begin{align}
        \begin{aligned}
            \timeDerivative{Z_V} =&\ - \kappa_2 \frac{(Z_V + K_2)\kappa_1 + \kappa_{-1}  - \sqrt{ \left((Z_V + K_2)\kappa_1 + \kappa_{-1}\right)^2 - 4 \kappa_1^2 Z_V K_2 } }{ 2\kappa_1 } \\
            \equiv&\ - \kappa_2 \eqb^{-}(K_2,Z_V)
            \end{aligned}
            \label{eq:FLLN_tQSSA_ODE}
        \end{align}
        with initial condition $Z_V(0)$ at time zero. In particular, if $Z_V(0)$ and $K_2$ are deterministic (non-random), then $Z_V$ is also deterministic, and hence as $n\rt \infty$, $\nZ_V \prt Z_V$, that is, 
        \begin{align*}
            \sup_{t\leq T}|\nZ_V(t)-Z_V(t)| \prt 0.
        \end{align*}
        \label{thm:tQSSA}
    \end{myTheorem}

\begin{myRemark}
   Since the path space of $Z_V$ is $C([0,T],\setOfPositiveReals)$, for deterministic  $Z_V$ the  convergence in probability in \Cref{thm:tQSSA} holds in the uniform metric by \cite[Chapter VI, Proposition 1.17]{jacod2003limit}.
\end{myRemark}
    
    \begin{proof}[Proof of \Cref{thm:tQSSA}]
   We divide the proof in three steps.\\
\np
{\em Step 1:}   By the hypothesis, \Cref{prop:rel_compactness}, and Prokhorov's theorem, the sequence $\{(\nK_2,\occ^{(n)}, Z_V^{(n)}) :  n \ge 1\}$ is relatively compact in $\setOfPositiveReals\times\spaceOfOccupationMeasures{T}{ \setOfPositiveReals\times  \setOfPositiveReals  \times [0, T]}\times D([0, T], \setOfPositiveReals)$. Let the triple $(K_2,\occ, Z_V)$ be a limit point of $\{(\nK_2,\occ^{(n)}, Z_V^{(n)}): n\ge 1 \}$, that is, there exists a subsequence --- which, by a slight abuse of notation, we also denote by $\{n\}$ --- along which 
   \begin{align}\label{eq:conv-occ-ZV-weak}
   (\nK_2,\occ^{(n)}, Z_V^{(n)}) \stackrel{n \rt \infty}\RT (K_2,\occ, Z_V).       
   \end{align}

   We will show for this limit point $Z_V$ is given by \eqref{eq:FLLN_tQSSA_ODE} and $\occ$  by 
   \begin{align}\label{eq:lim-occ-form-0}
   \occ(\differential{z_C}\times \differential{z_V}\times \differential{s})= \delta_{\eqb^{-}(K_2,Z_V(s))}(\differential{z_C})\delta_{Z_V(s)}(\differential{z_V})\differential{s},
   \end{align}
   that is, $\occ(A\times B \times [0,t]) = \int_0^t \delta_{\eqb^{-}(K_2,Z_V(s))}(A)\delta_{Z_V(s)}(B)\differential{s},$
   for any Borel sets $A,B \subset \setOfPositiveReals$.
    Since the differential equartion in \eqref{eq:FLLN_tQSSA_ODE} admits a unique solution on $[0,T]$, the limit point $(K_2,\occ, Z_V)$ is unique and independent of the subsequence.  Thus, the convergence holds along the entire sequence, which will complete the proof. 
We now work toward the above goal. 

\vs{.1cm}
\np{\em Step 2:} We now show that for any $g\in C^{(2)}_c(\setOfPositiveReals, \R),$ and $t \in [0,T]$,
\begin{align}
    \label{eq:occ-iden-0}
    \int_{\setOfPositiveReals\times  \setOfPositiveReals  \times [0, t]} \fgen_{K_2, z_V}g(z_C)\occ(\differential{z_C}\times \differential{z_V}\times \differential{s})=0, \quad \text{a.s.}
\end{align}

Fix  $g\in C^{(2)}_c(\setOfPositiveReals, \R)$. 
By the It\^o's formula for \acp{SDE} driven by \acp{PRM} (\eg, see \cite[Lemma~4.4.5]{Applebaum_2009Levy}, \cite[Theorem~5.1]{IkedaWatanabe2014Stochastic}) 
\begin{align}\label{eq:g-Z-C}
\begin{aligned}
    g(\nZ_C(t)) =&\ g(\nZ_C(0))+\int_{[0,\infty)\times [0, t]} \lf(g(\nZ_C(s-)+n^{-1}) -g(\nZ_C(s-))\ri)\\
    & \hs{.5cm}\times \indicator{[0, n^2 \kappa_1 (\nZ_V(s-)-\nZ_C(s-)) (\nK_2-\nZ_C(s-)]}{v} Q_1(\differential{v}\times \differential{s})\\
    &\ + \int_{[0,\infty) \times [0, t]} \lf(g(\nZ_C(s-)-n^{-1}) -g(\nZ_C(s-))\ri) \\
    &{}\hs{.5cm}\times
     \indicator{[0, n^2 \kappa_{-1}Z_C^{(n)} (s-)]}{v} Q_{-1}(\differential{v}\times \differential{s})\\
    &\ + \int_{[0,\infty) \times [0, t]} \lf(g(\nZ_C(s-)-n^{-1}) -g(\nZ_C(s-))\ri)\\
    & \hs{.5cm}\times
    \indicator{[0, n \kappa_{2}Z_C^{(n)} (s-)]}{v} Q_{2}(\differential{v}\times \differential{s})\\
    =&\ g(\nZ_C(0))+ n\int_0^t \fgen^{(n)}_{\nZ_V(s)} g(\nZ_c(s))\differential{s} +\mart^{(n)}_g(t)\\
    =&\ g(\nZ_C(0))+ n\int_{\R_{+}\times \R_{+}\times [0, t]} \fgen^{(n)}_{z_V} g(z_C)\occ^{(n)}(\differential{z_C}\times \differential{z_V}\times \differential{s})+\mart^{(n)}_g(t),
\end{aligned}
\end{align}
where for a fixed $z_V \in \setOfPositiveReals$ the operator $\fgen^{(n)}_{z_V}$ is defined by
\begin{align}
    \label{eq:fast-gen-n}
    \begin{aligned}
  \fgen^{(n)}_{z_V} g(z_C)=&\ n\kappa_1(z_V-z_C)(\nK_2-z_C)(g(z_C+n^{-1}) -g(z_C))\\
  &\ + (n\kappa_{-1}+\kappa_{2})z_C(g(z_C-n^{-1}) -g(z_C)),
   \end{aligned}
\end{align}
and the stochastic process $\mart^{(n)}_g$ is a zero-mean martingale given by
\begin{align}
    \label{eq:mart-g}
    \begin{aligned}
  \mart^{(n)}_g(t)=&\ \int_{\setOfPositiveReals \times [0, t]} \lf(g(\nZ_C(s-)+n^{-1}) -g(\nZ_C(s-))\ri)\\
    & \hs{.5cm}\times \indicator{[0, n^2 \kappa_1 (\nZ_V(s-)-\nZ_C(s-)) (\nK_2-\nZ_C(s-))]}{v} \tilde Q_1(\differential{v}\times \differential{s})\\
    &\ + \int_{\setOfPositiveReals \times [0, t]} \lf(g(\nZ_C(s-)-n^{-1}) -g(\nZ_C(s-))\ri) \\
    &\hs{.5cm}\times
    \indicator{[0, n^2 \kappa_{-1}Z_C^{(n)} (s-)]}{v} \tilde Q_{-1}(\differential{v}\times \differential{s})\\
    &\ + \int_{\setOfPositiveReals\times [0, t]} \lf(g(\nZ_C(s-)-n^{-1}) -g(\nZ_C(s-))\ri) \\
    &\hs{.5cm}\times
    \indicator{[0, n \kappa_{2}Z_C^{(n)} (s-)]}{v} \tilde Q_{2}(\differential{v}\times \differential{s}). 
    \end{aligned}
\end{align}
Here $\tilde Q_1, \tilde Q_{-1}, \tilde Q_2$ are the compensated \acp{PRM} corresponding to $Q_1, Q_{-1}, Q_2$ respectively.
Write 
\begin{align*}
\fgen^{(n)}_{z_V}g(z_C) \equiv \fgen_{\nK_2, z_V}g(z_C)+ \err^{(n)}_{g,1}(z_C,z_V),
\end{align*}
where recall that for $\nK_2>0, z_V \in \setOfPositiveReals$, the operator $\fgen_{\nK_2, z_V} $ is given by \eqref{eq:fast_gen}. Notice that  by second-order Taylor expansion and the fact that $\partial^2 g$ is bounded (since $g \in C^{(2)}_c(\setOfPositiveReals, \R)$), the error term $\err^{(n)}_{g,1}(z_C,z_V) \equiv\fgen^{(n)}_{z_V}g(z_C) -\fgen_{\nK_2, z_V}g(z_C) $ can be estimated as
\begin{align*}
 |\err^{(n)}_{g,1}(z_C,z_V)| \leq&\  \const_{g,0}\lf( n^{-1}\kappa_1(z_V-z_C)(\nK_2-z_C) 
  + (\kappa_{-1}n^{-1}+ \kappa_2n^{-1}+\kappa_2n^{-2})z_C\ri)\\
 \leq &\ \const_{g,1}n^{-1}(1+\nK_2)(1+z^2_V+z_C^2),
\end{align*}
for suitable constants, $\const_{g,0}$, and  $\const_{g,1}$. Therefore, by the tightness of the sequence of random variables $\{\nK_1: n\ge 1 \}$, as $n \rt \infty$, we have 
\begin{align*}
   \int_{\R_{+} \times \R_{+} \times [0, T]}& |\err^{(n)}_{g,1}(z_C,z_V)| \occ^{(n)}(\differential{z_C}\times \differential{z_V}\times \differential{s}) \\
   \leq&\ \const_{g,1}(1+\nK_2)n^{-1}\int_0^T (1+(\nZ_V(s))^2+(\nZ_C(s))^2)\ \differential{s}\\
   \leq&\ \const_{g,1}(1+\nK_2)(1+2(\nK_1)^2)Tn^{-1} \prt 0.
\end{align*}
Next, by the Lipschitz continuity of $g$ and the \ac{BDG} inequality  \citep[Appendix D]{Budhiraja2019Analysis}, we can estimate the predictable quadratic variation of the martingale $\mart^{(n)}_g$ as
\begin{align*}
\<\mart^{(n)}_g\>_T \leq \const_{g,2}(1+(\nZ_V(T))^2+(\nZ_C(T))^2) \leq 2\const_{g,2}(1+(\nK_1)^2)\eqcomma 
\end{align*}
for some constant $\const_{g,2}>0$. Thus, $\<\mart^{(n)}_g\>_T/n^2 \prt 0$ as $n\rt \infty$, and hence, as in the proof of \Cref{prop:rel_compactness} by the Lenglart--Rebolledo inequality, we have as $n\rt \infty$ 
\begin{align*}
    \sup_{t\leq T}|\mart^{(n)}_g(t)|/n \prt 0.
\end{align*}
Now from the decomposition in \eqref{eq:g-Z-C} we have
\begin{align}\label{eq:conv-occ-fgen-0}
\begin{aligned}
    &\int_{\R_{+} \times \R_{+} \times [0, t]} \fgen_{\nK_2, z_V}g(z_C)\occ^{(n)}(\differential{z_C}\times \differential{z_V}\times \differential{s})\\
    & = n^{-1}\Big(g(\nZ_C(t))-g(\nZ_C(0))-\mart^{(n)}_g(t)\Big) \\
    &\quad 
    - \int_{\R_{+} \times \R_{+} \times [0, t]} \err^{(n)}_{g,1}(z_C,z_V) \occ^{(n)}(\differential{z_C}\times \differential{z_V}\times \differential{s}) 
    \prt 0, \quad \mbox{as } n \rt \infty.
\end{aligned}
\end{align}
Next, notice that the mapping $(K, z_V, z_C) \rt \fgen_{K, z_V}g(z_C)$ satisfies
\begin{align*}
|\fgen_{K, z_V}g(z_C)\ -\fgen_{K', z_V}g(z_C)| \leq \kappa_1\|\partial g\|_{\infty}|K-K'|(z_V+z_C),
\end{align*}
and hence, satisfies the hypothesis (on $h$) in \Cref{lem:conv-occ-ZV}:\eqref{item:int-occ-conv}. By virtue of \eqref{eq:conv-occ-ZV-weak} and \Cref{lem:conv-occ-ZV}:\eqref{item:int-occ-conv}, 
\begin{align} \label{eq:conv-occ-fgen-1}
    \begin{aligned}
    & \int_{\R_{+} \times \R_{+} \times [0, t]} \fgen_{\nK_2, z_V}g(z_C)\occ^{(n)}(\differential{z_C}\times \differential{z_V}\times \differential{s}) \\
    &\quad \quad\quad  \nRT \int_{\R_{+}\times \R_{+}\times [0, t] } \fgen_{K_2, z_V}g(z_C)\occ(\differential{z_C}\times \differential{z_V}\times \differential{s}),
    \end{aligned}
\end{align}
which, along with \eqref{eq:conv-occ-fgen-0}, establishes \eqref{eq:occ-iden-0}.

\vs{.1cm}
\np{\em Step 3:} We now identify the measure $\occ$. Notice that because of \eqref{eq:supp-occ-n},
\begin{align}\label{eq:supp-occ-lim}
    \mathrm{supp}(\occ) \subset \{(z_C,z_V) \in \setOfPositiveReals\times \setOfPositiveReals: z_C \leq z_V \}\times \setOfPositiveReals,
    \end{align}
Next, let $\occ_{(2,3)}$ and $\occ^{(n)}_{(2,3)}$ be the marginal distributions of $z_V$ and the time component corresponding to $\occ$ and $\occ^{(n)}$, respectively; that is, define $\occ_{(2,3)}$ by  $\occ_{(2,3)}(B\times [0,t])\equiv\occ(\setOfPositiveReals\times B\times [0,t])$ ($\occ^{(n)}_{(2,3)}$ is similarly defined).
Now, decompose the measure $\occ$ as 
  \begin{align}\label{eq:occ-decom}
  \occ(\differential{z_C}\times \differential{z_V}\times \differential{s}) = \occ_{(1|2,3)}(\differential{z_C}|z_V,s) \occ_{(2,3)}(\differential{z_V}\times \differential{s})\eqstop 
  \end{align}
 Note that for each fixed $z_V\geq 0$ and $s\geq 0$,
\begin{align}\label{eq:supp-occ-C-lim}
    \mathrm{supp}(\occ_{(1|2,3)}(\cdot|z_V,s)) \subset [0,z_V].
    \end{align}
 Clearly, the convergence in \eqref{eq:conv-occ-ZV-weak} implies that as $n\rt \infty$, $\occ^{(n)}_{(2,3)} \RT \occ_{(2,3)}$  in $\spaceOfOccupationMeasures{T}{\setOfPositiveReals\times [0, T]}.$ We first argue that 
\begin{align}
    \label{eq:occ-ZV-form}
  \occ_{(2,3)}(\differential{z_V}\times \differential{s}) \stackrel{d}= \delta_{Z_V(s)}(\differential{z_V})\differential{s}.
  \end{align}
Fix any $h \in C_b(\setOfPositiveReals\times [0,T],\R)$. By (a simpler version) of \Cref{lem:conv-occ-ZV}:\eqref{item:int-occ-conv}  (or just by the continuous mapping theorem)
$$\int h\ \differential{\occ^{(n)}_{(2,3)}} \stackrel{n\rt \infty} \RT \int_{\setOfPositiveReals \times [0, t]} h\ \differential{ \occ_{(2,3)}},$$
where for $\nu \in \spaceOfOccupationMeasures{T}{\setOfPositiveReals \times [0,T]}$, $\int h \differential{\nu} \equiv \int_{\setOfPositiveReals \times [0, t]} h(z_V,s) \nu(\differential{z_V}\times \differential{s}).$

On the other hand, observing that the function $\phi \in D([0,T],\setOfPositiveReals) \rt \int_0^T h(\phi(s),s)\differential{s}$ is continuous \footnote{$\phi_n \rt \phi$ in Skorohod topology implies $\phi_n(s) \rt \phi(s)$ for all continuity points of $\phi$. Since $\phi$ is c\`adl\`ag, the set of discontinuity points of $\phi$ is  countable. Hence, by the dominated convergence theorem $\int_0^T h(\phi_n(s),s)ds \rt \int_0^T h(\phi(s),s)ds$. }, we have by the continuous mapping theorem,
\begin{align*}
 \int h\ \differential{\occ^{(n)}_{(2,3)}}  =&\ \int_0^T h(\nZ_V(s),s)\differential{s}   \stackrel{n\rt \infty} \RT  \ \int_0^t h(Z_V(s),s)\differential{s}\\
    =&\ \int_{[0,t]\times \setOfPositiveReals} h(z_V) \delta_{Z_V(s)}(\differential{z_V})\differential{s}.
\end{align*}
Hence,
\begin{align*}
\int h(z_V,s) \occ_{(2,3)}(\differential{z_V}\times \differential{s}) \stackrel{d} = \int h(z_V,s) \delta_{Z_V(s)}(\differential{z_V})\differential{s}, \quad h \in C_b(\setOfPositiveReals \times [0,T],\R).
\end{align*}
It readily follows by the Cram\'er--Wold theorem that  for any finite collection $\{h_1,h_2,\hdots,h_m\} \subset C_b(\setOfPositiveReals \times [0,T],\R),$
\begin{align}\label{eq:int-dist-equality}
\begin{aligned}
\Big(& \int_{\setOfPositiveReals \times [0, T]} h_j(z_V,s) \occ_{(2,3)}(\differential{z_V}\times \differential{s});\ j=1\hdots m\Big) \\
& \stackrel{d} = \Big(\int_{\setOfPositiveReals\times [0, T]} h_j(z_V,s) \delta_{Z_V(s)}(\differential{z_V})\differential{s};\ j=1\hdots m\Big)
\end{aligned}
\end{align}
Since $\mathcal{B}(\spaceOfOccupationMeasures{T}{\setOfPositiveReals\times [0, T]})$, the Borel $\s$-field on $\spaceOfOccupationMeasures{T}{\setOfPositiveReals\times [0, T]}$, is generated by the sets of the form 
$$\lf\{\nu: \int h_j \differential{\nu} < a_j,\ j=1,2\hdots, m \ri\}, \quad a_j \in \R,\ h_j \in C_b(\setOfPositiveReals \times [0,T],\R),$$
and the collection of such sets form a $\pi$-system, we
conclude that \eqref{eq:int-dist-equality} implies \eqref{eq:occ-ZV-form}. Since convergence in distribution determines the limit only up to equality in distribution, and we have identified the distribution of limiting $\occ_{(2,3)}$ in \eqref{eq:occ-ZV-form}, we  write $\occ^{(n)}_{(2,3)} \RT \occ_{(2,3)}$, with
\begin{align}
    \label{eq:occ-ZV-form-2}
  \occ_{(2,3)}(\differential{z_V}\times \differential{s}) = \delta_{Z_V(s)}(\differential{z_V})\differential{s}.
\end{align}

We next identify the measure $\occ_{(1|2,3)}(\cdot|z_V,s)$. 
By the decomposition in \eqref{eq:occ-decom} and \eqref{eq:occ-ZV-form-2}, we have 
\begin{align*}
\int_{[0, t]\times \R_{+} \times \R_{+} }& \fgen_{K_2, z_V}g(z_C)\occ_{(1|2,3)}(\differential{z_C}|z_V,s)\delta_{Z_V(s)}(\differential{z_V})\differential{s}\\
=&\ \int_{[0, t]\times \R_{+} \times \R_{+}}\fgen_{K_2, Z_V(s)}g(z_C)\occ_{(1|2,3)}(\differential{z_C}|Z_V(s),s)\differential{s}\ =\ 0.
\end{align*}
Since this is true for all $t \in [0,T]$ and the space $C^{(2)}_c(\setOfPositiveReals, \R)$ is separable, we have for a.a. $s \in [0,T]$ that 
\begin{align}\label{eq:conv-occ-fgen-2}
    \int_{\R_{+}}\fgen_{K_2, Z_V(s)}g(z_C)\occ_{(1|2,3)}(\differential{z_C}|Z_V(s),s)=0 ,\ \mbox{ for any } g \in C^{(2)}_c(\setOfPositiveReals, \R).
\end{align}
Because of \eqref{eq:supp-occ-C-lim},  we can apply \Cref{lem:fast-inv-dist} to conclude that for a.a. $s \in [0,T]$,  $\occ_{(1|2,3)}(\cdot|Z_V(s),s) = \delta_{\eqb^{-}(K_2,Z_V(s))}$. Consequently, by \eqref{eq:occ-decom}, \eqref{eq:occ-ZV-form-2}, we see that for the limit point $(Z_V, \occ)$, \eqref{eq:lim-occ-form-0} holds.
Now,  recall that \eqref{eq:ZV-split} gives
\begin{align*}
   \nZ_V(t) - \nZ_V(0) + \int_{0}^{t} \kappa_2 \nZ_C(s) \differential{s} = \mart^{(n)}_V(t).
\end{align*}
Because of \eqref{eq:conv-occ-ZV-weak}, \Cref{lem:conv-occ-ZV}, \eqref{eq:lim-occ-form-0}, and the continuous mapping theorem show that as $n\rt \infty$,
\begin{align*}
    \nZ_V(t)\ - &\ \nZ_V(0) + \int_{0}^{t} \kappa_2 \nZ_C(s) \differential{s} \\
    & \RT  Z_V(t) - Z_V(0) + \int_{\R_{+}\times \R_{+} \times [0, t] } \kappa_2 z_C \occ(\differential{z_C}\times \differential{z_V}\times \differential{s})\\
    & = Z_V(t) - Z_V(0)+\int_0^t \kappa_2 \eqb^{-}(K_2,Z_V(s))  \differential{s}.
\end{align*} 
On the other hand, from the proof of \Cref{prop:rel_compactness},  $\sup_{t\leq T}|\mart^{(n)}_V(t)| \prt 0$, as $n\rt \infty$ (see \eqref{eq:mart-conv}). It follows that the limit point $Z_V$ must satisfy the \ac{ODE} \eqref{eq:FLLN_tQSSA_ODE}, which establishes the goal outlined in  {\em Step 1}.

\end{proof}

    \begin{figure}[t]
        \centering
        \includegraphics[width=0.75\textwidth]{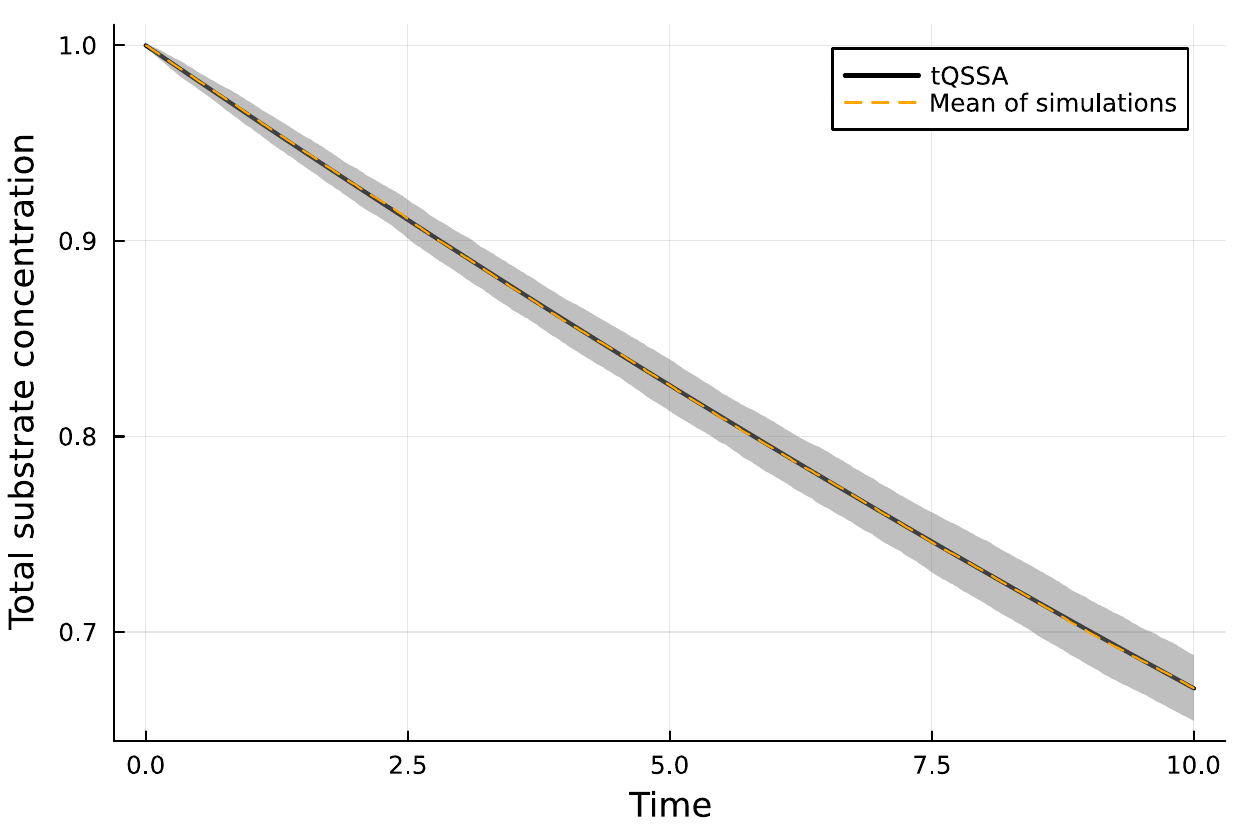}
        \caption{Accuracy of the \ac{tQSSA} for the \ac{MM} enzyme kinetic model. The plot shows that the exact Doob--Gillespie simulations \citep{Wilkinson2018SMS} and the \ac{tQSSA} solution for the \ac{MM} enzyme kinetic model are very close. The parameters used in the simulations are $\kappa_1 = 1$, $\kappa_{-1} = 1$, $\kappa_2 = 0.75$, $K_2 = 0.1$, $K_1=1.0$, and $n=1000$.}
        \label{fig:tQSSA}
    \end{figure}

    See \Cref{fig:tQSSA} for a comparison of exact Doob--Gillespie simulations (see \cite{Wilkinson2018SMS}) with the \ac{tQSSA} solution for the \ac{MM} enzyme kinetic model. 

    \section{The \acl{FCLT}}
    \label{sec:FCLT}
  We now study the fluctuations around the \ac{tQSSA}. Toward this end,  define the scaled stochastic process 
    \begin{align*}
        \nU_V(t) \defeq \sqrt{n}(\nZ_V(t) - Z_V(t))\eqcomma 
    \end{align*}
    where $Z_V$ is the limit point of $\nZ_V$ solving the \ac{ODE} in \eqref{eq:FLLN_tQSSA_ODE} as established in \Cref{thm:tQSSA}. The proof of the \ac{FCLT} hinges on the solution of a Poisson equation, which we now study.

    \subsection{Solution of a Poisson equation and its properties}

    For each $z_V \in \setOfPositiveReals$, let $F(z_V, \cdot): \setOfPositiveReals \rt \R$ be a solution of the Poisson equation 
    \begin{align}
        \label{eq:poi}
        \fgen_{K_2,z_V}F(z_V,\cdot)(z_C) = - (z_C - \eqb^-(K_2,z_V)), \quad \text{ for } z_C \leq z_V
    \end{align}
    under the centering condition: $F(z_V,\eqb^-(K_2,z_V))=0.$ Using the fact that  $\eqb^\pm(K_2,z_V)$ are the two roots of  \eqref{eq:ode-eqb-eqn}, we see that $\fgen_{K_2,z_V}F(z_V,\cdot)(z_C)$ can be factorized as
    \begin{align*}
    \fgen_{K_2,z_V}F(z_V,\cdot)(z_C) = \kappa_1(z_C -\eqb^-(K_2,z_V))(z_C -\eqb^+(K_2,z_V))\partial_2F(z_V,z_C),    \end{align*}
    whence it readily follows that for $z_C \leq z_V$, the solution $F$ is given by
   \begin{align}
    \begin{aligned}
        F(z_V,z_C)=&\ -\f{1}{\kappa_1}\ln\lf(\f{\eqb^+(K_2,z_V)-z_C}{\eqb^+(K_2,z_V)-\eqb^-(K_2,z_V)}\ri),\\
    =&\ -\f{1}{\kappa_1} \lf[\ln(\eqb^+(K_2,z_V)-z_C) - \f{1}{2}\ln(D(K_2,z_V)) + \ln(\kappa_1)\ri],
    \end{aligned}
    \label{eq:Poisson_soln}
\end{align}
 where the discriminant $D(K_2,z_V)$ is defined in \eqref{eq:disc}.
Clearly, from the expression of $\eqb^+(K_2,z_V)$ and the observations in  \eqref{eq:min-dist-eqb} and \eqref{eq:disc}, we have the following upper estimate for the function $F$:
\begin{align}
    \label{eq:poi-growth}
    |F(z_V,z_C)| \leq \const_{F,0}\lf(1+\ln(1+z_V+z_C)\ri)
\end{align}
for some constant $\const_{F,0}>0.$
Moreover, notice that the derivatives of $F$ with respect to $z_C$ can be estimated as follows: for $z_C\leq z_V$
  \begin{align}
        \label{eq:poi-sol-2nd partial-1}
        \partial_2 F(z_V,z_C)=&\ \f{1}{\kappa_1(\eqb^+(K_2,z_V)-z_C)} \leq 2/\kappa_{-1},\\
        \label{eq:poi-sol-2nd partial-2}
        \partial^2_2 F(z_V,z_C)=&\ \f{1}{\kappa_1(\eqb^+(K_2,z_V)-z_C)^2} \leq 4\kappa_1/\kappa^2_{-1}.
    \end{align}
   Also, note that the derivative of $F$ with respect to $z_V$ is given by 
   \begin{align*}      
      \partial_1 F(z_V,z_C)=&\ -\f{1}{\kappa_1}\left[(\eqb^+(K_2,z_V)-z_C)^{-1}\partial \eqb^+(K_2,z_V) \right. \\
      &\quad \left. 
        - (2D(K_2,z_V))^{-1}\partial D(K_2,z_V)\right] \\
      =&\ -\f{1}{\kappa_1}\big[(\eqb^+(K_2,z_V)-z_C)^{-1}\partial \eqb^+(K_2,z_V)\\
      &\ \hs{1cm}- (D(K_2,z_V))^{-1}((z_V-K_2)\kappa_1^2+\kappa_1\kappa_{-1})\big].
   \end{align*}
   Therefore,  we have 
   \begin{align}
   \label{eq:poi-sol-1st partial-1}     
   |\partial_1 F(z_V,z_C)| \leq \const_{F,1}(1+z_V).
\end{align}
Next,   notice that  $\gen^{(n)}$ defined in \eqref{eq:reduced_Z_generator} admits the decomposition:
    \begin{align}
    \label{eq:gen-expansion}
    \begin{aligned}
     \gen^{(n)} F(z_V, z_C)=&\ n \lf(\fgen^{(n)}_{z_V}F(z_V,\cdot)(z_C)+\err^{(n)}_{F,0}(z_V, z_C)\ri)\\
     =&\ n\lf(\fgen_{\nK_2, z_V}F(z_V,\cdot)(z_C)+\err^{(n)}_{F,0}(z_V,z_C)+\err^{(n)}_{F,1}(z_V, z_C)\ri)\\
     =&\ n\left(\fgen_{K_2, z_V}F(z_V,\cdot)(z_C)+\err^{(n)}_{F,0}(z_V,z_C)+\err^{(n)}_{F,1}(z_V, z_C) \right. \\
     &\quad \left. +\err^{(n)}_{F,2}(z_V, z_C)\right),
     \end{aligned}
    \end{align}
    where
    \begin{align*}
      \err^{(n)}_{F,0}(z_V, z_C)=&\ \kappa_2z_C\lf(F(z_V-n^{-1}, z_C-n^{-1})-F(z_V, z_C-n^{-1})\ri)\eqcomma  \\
      \err^{(n)}_{F,1}(z_V, z_C) \equiv&\ \fgen^{(n)}_{z_V}F(z_V,\cdot)(z_C) - \fgen_{\nK_2, z_V}F(z_V,\cdot)(z_C)\eqcomma \\
      \err^{(n)}_{F,2}(z_V, z_C) \equiv&\ \fgen_{\nK_2, z_V}F(z_V,\cdot)(z_C)-\fgen_{K_2, z_V}F(z_V,\cdot)(z_C)\\
            =&\ \kappa_1(z_V-z_C)(\nK_2-K_2)\partial_2F(z_V,z_C)\eqstop 
    \end{align*}
    Note that Taylor expansion, \eqref{eq:poi-sol-1st partial-1}, \eqref{eq:poi-sol-2nd partial-2} and \eqref{eq:poi-sol-2nd partial-1} show for some constants $\const_{F,0}, \const_{F,1}>0$ that 
    \begin{align}\label{eq:err-bds}
      \begin{aligned}
    | \err^{(n)}_{F,0}(z_V, z_C)| \leq&\ \const_{F,0}(1+z_V)z_Cn^{-1}\eqcomma \\
    | \err^{(n)}_{F,1}(z_V, z_C)| \leq&\ \big|n^{-1}\kappa_1(z_V-z_C)(\nK_2-z_C)+(n^{-1}\kappa_{-1}+\kappa_2n^{-2})z_C\big| \\
    &\quad \times  
     \big|\sup_{z_C \leq z_V}\partial_2^2F(z_V,z_C)\big|\\
    \leq&\ \const_{F,1}(1+\nK_2)(1+z_V+z_C)^2n^{-1} \eqcomma \\
     | \err^{(n)}_{F,2}(z_V, z_C)| \leq&\ 2\kappa_1(\kappa_{-1})^{-1}|\nK_2-K_2|(z_V+z_C).
      \end{aligned}
    \end{align}

    We will now present the \ac{FCLT} for the process $\nU_V$.
    \begin{myTheorem}
        Suppose that as $n\rt \infty$,  
        $$(\nZ_V(0), \nU_V(0), \nK_2,\sqrt{n} (\nK_2 -K_2))  \RT (Z_V(0), U_V(0),K_2,  \tilde K_2),$$ 
        where $Z_V(0)$ and $K_2$ are deterministic, but $(U_V(0), \tilde K_2)$ can be random. Also assume that $\sup_n\EE(\nK_1)<\infty$. Then, the sequence  $\nU_V$ is $C$-tight in $D([0,T], \R)$, and as $n\rt \infty$, 
        $$\nU_V \RT U_V,$$
        where $U_V$ satisfies the \ac{SDE}
        \begin{align}\label{eq:FCLT}
       \non
            U_V(t) =&\ U_V(0) - \kappa_2  \tilde K_2\int_0^t \f{Z_V(s)-\eqb^-(K_2,Z_V(s))}{\sqrt{ \left((Z_V(s) + K_2)\kappa_1 + \kappa_{-1}\right)^2 - 4 \kappa_1^2K_2Z_V(s)}} \differential{s}\\
            &\ +\kappa_2\int_0^t \partial\eqb^{-}(K_2,Z_V(s))U_V(s)\differential{s}
             + \int_{0}^{t}\sqrt{{\kappa_2}\eqb^{-}(K_2,Z_V(s))}\differential{W}(s),
        \end{align}
         the stochastic process $W$ is the standard Brownian motion, and $\eqb^-(K_2,z_V)$ and $\partial\eqb^{-}(K_2,z_V)$ are as in \eqref{eq:ode-eqb-pts} and \eqref{eq:deriv-eqb}, respectively.
        \label{thm:Z_V_FCLT}
    \end{myTheorem}
    \begin{proof}[Proof of \Cref{thm:Z_V_FCLT}]
    First note that \eqref{eq:ZV-split} and \eqref{eq:FLLN_tQSSA_ODE} give
    \begin{align}
    \label{eq:U-exp-0}
    \begin{aligned}
      \nU_V(t) =&\ \sqrt{n}(\nZ_V(t)-Z_V(t)) \\
       =&\ \nU_V(0)-\kappa_2\sqrt n\int_0^t \lf(\nZ_C(s) - \eqb^{-}(\nZ_V(s))\ri)\differential{s}\\
      &\ +\kappa_2\sqrt n\int_0^t \lf(\eqb^{-}(\nZ_V(s)) - \eqb^{-}(Z_V(s))\ri) \differential{s} + \sqrt{n} \mart^{(n)}_V(t).
    \end{aligned}
    \end{align}
    
  \np 
We first consider the {\em second term} on the right-side of \eqref{eq:U-exp-0}. Applying the It\^{o}'s formula \citep[Theorem~5.1]{IkedaWatanabe2014Stochastic} with  $F$ as the solution of the Poisson equation \ref{eq:poi}, we have 
    \begin{align*}
        F(\nZ_V(t), \nZ_C(t)) - F(\nZ_V(0), \nZ_C(0)) = \nM_F(t) + \int_{0}^{t} \gen^{(n)} F(\nZ_V(u), \nZ_C(u))\differential{u}\eqcomma 
    \end{align*}
    where the zero-mean martingale $\nM_F$ is given by
    \begin{align*}
        \nM_F(t) \defeq& \int_{[0,\infty) \times [0, t]} \left(F(\nZ_V(s-)-n^{-1}, \nZ_C(s-)-n^{-1}) - F(\nZ_V(s-), \nZ_C(s-)) \right) \\
        &{}\quad\quad  \times 
        \indicator{[0, n\kappa_2 \nZ_C(s-)]}{v}\tilde{Q}_2(\differential{v}\times \differential{s})\\
        &\quad + \int_{[0,\infty) \times [0, t]} \left(F(\nZ_V(s-), \nZ_C(s-)+n^{-1}) - F(\nZ_V(s-), \nZ_C(s-)) \right)\\
        &\quad\quad  \times \indicator{[0, n^2\kappa_1 (\nZ_V(s-) - \nZ_C(s-))(\nM -\nZ_C(s-))]}{v}\tilde{Q}_1(\differential{v}\times \differential{s})\\
        &\quad + \int_{[0,\infty) \times [0, t]} \left(F(\nZ_V(s-), \nZ_C(s-)- n^{-1}) - F(\nZ_V(s-), \nZ_C(s-)) \right) \\
        &{}\quad \quad \times 
        \indicator{[0, n^2\kappa_{-1}\nZ_C(s-) ]}{v}\tilde{Q}_{-1}(\differential{v}\times \differential{s}),
    \end{align*}
    and $\tilde{Q}_1, \tilde{Q}_{-1}, \tilde{Q}_2$ are the compensated \acp{PRM} associated with $Q_1, Q_2$ and $Q_{-1}$. By \eqref{eq:poi} and \eqref{eq:gen-expansion}, we have 
    \begin{align*}
      \int_{0}^{t} \gen^{(n)} F(\nZ_V(s), \nZ_C(s))\differential{s}\ =&\ -n \int_{0}^{t} \left( \nZ_C(s) - \eqb^-(K_2,\nZ_V(s)) \right)\differential{s} \\
        &\ + n\sum_{i=0}^2\int_{0}^{t} \err^{(n)}_{F,i}\lf(\nZ_V(s), \nZ_C(s)\ri)\differential{s}\eqstop
    \end{align*}
    Thus, it follows that 
    \begin{align}\label{eq:poieq-expnsn}
    \begin{aligned}
     -n^{1/2} \int_{0}^{t}& \left( \nZ_C(s) - \eqb^-(K_2,\nZ_V(s)) \right)\differential{s}\\
     =&\ n^{-1/2}\lf(F(\nZ_V(t), \nZ_C(t)) - F(\nZ_V(0), \nZ_C(0))\ri) -n^{-1/2}\nM_F(t)\\
     &\ -n^{1/2} \sum_{i=0}^2\int_{0}^{t} \err^{(n)}_{F,i}\lf(\nZ_V(s), \nZ_C(s)\ri)\differential{s}.
    \end{aligned}
    \end{align}
    Notice that by virtue of the estimates in \eqref{eq:err-bds} and the conservation laws in \eqref{eq:conservation_laws} the following estimates hold:
  \begin{align}\label{eq:var-err-est}
   \begin{aligned}
      \sup_{s\leq T} |\err^{(n)}_{F,0}(\nZ_V(s), \nZ_C(s))|\ \leq&\ \const_{F,0}(1+\nK_1)\nK_1n^{-1}\eqcomma \\
      \sup_{s\leq T} |\err^{(n)}_{F,1}(\nZ_V(s), \nZ_C(s))|\ \leq&\ \const_{F,2}(1+\nK_2)(1+\nK_1)^2n^{-1}\eqcomma \\
     \sup_{s\leq T} |\err^{(n)}_{F,2}(\nZ_V(s), \nZ_C(s))|\ \leq&\  4\kappa_1(\kappa_{-1})^{-1}\nK_1|\nK_2-K_2|.
  \end{aligned}
  \end{align}
  Also notice that by \eqref{eq:poi-growth} and \eqref{eq:conservation_laws}, we have 
\begin{align*}
    \sup_{t\leq T}|F(\nZ_V(t), \nZ_C(t)) - F(\nZ_V(0), \nZ_C(0))| \leq & 2\const_{F,0}\lf(1+\ln(1+2\nK_1)\ri)\eqstop 
\end{align*}
Therefore, the following limits hold as $n\rt \infty$
\begin{align}
\label{eq:FE-conv-0}
  \begin{aligned}
 n^{-1/2}\sup_{t\leq T}|F(\nZ_V(t), \nZ_C(t)) - F(\nZ_V(0), \nZ_C(0))| \prt& \ 0\eqcomma \\
 \sup_{s\leq T} n^{1/2} |\err^{(n)}_{F,i}(\nZ_V(s), \nZ_C(s))| \prt&\ 0, \quad i=0,1\\ 
 n^{1/2}\int_0^T |\err^{(n)}_{F,i}(\nZ_V(s), \nZ_C(s))| \differential{s} \prt& \ 0, \quad i=0,1.\\
  \end{aligned}
\end{align}
The tightness of the processes $n^{1/2}\int_0^\cdot \err^{(n)}_{F,2}(\nZ_V(s), \nZ_C(s))\differential{s}$ in $C([0,T],\R)$ follows from the last inequality of \eqref{eq:var-err-est}, the assumption of tightness of the sequence $n^{1/2}(\nK_2-K)$ and $\{\nK_1\}$ (since $\sup_n\EE(\nK_1)<\infty$), and \Cref{lem:int-proc-tight} in Appendix~\ref{sec:aux-results}.

We next consider the {\em third term} on the right-side of \eqref{eq:U-exp-0}. Applying first-order Taylor expansion (mean value theorem) to the integrand of the third term on right side of \eqref{eq:U-exp-0} we can write 
 \begin{align}
     \label{eq:eqb-expansion-1}
     \begin{aligned}
        \int_0^t \lf(\eqb^{-}(\nZ_V(s)) - \eqb^{-}(Z_V(s))\ri) \differential{s} =&\ n^{-1/2} \int_0^t\SC{D}_{C,*}(\nZ_V(s),Z_V(s))\nU_V(s)\ \differential{s}\eqcomma 
    \end{aligned}
 \end{align}
 where, because of \eqref{eq:eqb-deriv-lin}, the term $\SC{D}_{C,*}(\nZ_V(s),Z_V(s))$ can be estimated as
 \begin{align*}
    |\SC{D}_{C,*}(\nZ_V(s),Z_V(s))| \leq \const_{*,1}(1+ \nZ_V(s)+Z_V(s))\eqstop 
 \end{align*}
 By the conservation laws in \eqref{eq:conservation_laws} and the continuity of the function $Z_V$, we have the following estimate  
  \begin{align*}
    \sup_{s\leq T}|\SC{D}_{C,*}(\nZ_V(s),Z_V(s))| \leq \const_{*,1}\lf(1+ \nK_1+\sup_{s\leq T}Z_V(s)\ri).
 \end{align*}  
 Therefore, by the hypothesis, the collection $\{\sup_{s\leq T}|\SC{D}_{C,*}(\nZ_V(s),Z_V(s))|\}$ is a tight sequence of $\setOfPositiveReals$-valued random variables. Plugging \eqref{eq:poieq-expnsn} and \eqref{eq:eqb-expansion-1}  in the expression of $\nU_V$ in  \eqref{eq:U-exp-0}, we get 
    \begin{align}\label{eq:Un-SDE-0}
    \begin{aligned}
        \nU_V(t) & = \nU_V(0) + n^{-1/2}\kappa_2\lf(F(\nZ_V(t), \nZ_C(t)) - F(\nZ_V(0), \nZ_C(0))\ri)\\
        &\quad  +\kappa_2\int_0^t\SC{D}_{C,*}(\nZ_V(s),Z_V(s))\nU_V(s) + \mar^{(n)}(t)\\
        &\quad -n^{1/2}\kappa_2 \sum_{i=0}^2\int_{0}^{t} \err^{(n)}_{F,i}\lf(\nZ_V(s), \nZ_C(s)\ri)\differential{s} \eqcomma 
    \end{aligned}    
    \end{align} 
    where $\mar^{(n)}$
    is a zero-mean martingale given by
    \begin{align*}
        \mar^{(n)}(t) &{}\equiv\ n^{1/2}\mart^{(n)}_V(t)-n^{-1/2}\kappa_2\nM_F(t)  \\
        &{} =  -\frac{1}{\sqrt{n}} \Big[ \int_{\setOfPositiveReals \times [0, t]} \big(1 + \kappa_2\left(F(\nZ_V(s-)-n^{-1}, \nZ_C(s-)-n^{-1}) \right. \\
        &{}\quad\quad\quad\quad  \left. 
        - F(\nZ_V(s-), \nZ_C(s-))\big)\right) 
        \indicator{[0, n\kappa_2 \nZ_C(s-)]}{v}\tilde{Q}_2(\differential{v}\times \differential{s}) \\
        &\quad - \kappa_2\int_{\setOfPositiveReals\times [0, t]} \left(F(\nZ_V(s-), \nZ_C(s-)-n^{-1}) - F(\nZ_V(s-), \nZ_C(s-)) \right)\\
        &\quad\quad \quad \times \indicator{[0, n^2\kappa_1 (\nZ_V(s-) - \nZ_C(s-))(\nK_2 -\nZ_C(s-))]}{v}\tilde{Q}_1(\differential{v}\times \differential{s})\\
        &\quad  - \kappa_2\int_{\setOfPositiveReals\times [0, t]} \left(F(\nZ_V(s-), \nZ_C(s-)- n^{-1}) - F(\nZ_V(s-), \nZ_C(s-)) \right) \\
        &{}\quad \quad \quad 
         \times \indicator{[0, n^2\kappa_{-1}\nZ_C(s-) ]}{v}\tilde{Q}_{-1}(\differential{v}\times \differential{s})\Big]\eqstop \\
    \end{align*}

   Note that by the first-order Taylor expansion, the predictable quadratic variation of $\mar^{(n)}$ is given by
   
   \begin{align}\label{eq:pred-quad-mart-F}
   \begin{aligned}
        \predictableVariation{\mar^{(n)}}_t &= \int_{0}^{t} \left(1 + \kappa_2\left(F(\nZ_V(s)-n^{-1}, \nZ_C(s)-n^{-1}) - F(\nZ_V(s), \nZ_C(s))\right)\right)^2 \\
        &\quad\quad \times \kappa_2 \nZ_C(s)\differential{s}\\
        &\quad + \kappa_2^2\int_{0}^{t} n\left(F(\nZ_V(s), \nZ_C(s)-n^{-1}) - F(\nZ_V(s), \nZ_C(s)) \right)^2 \\
        &{}\quad\quad\quad 
        \times  \kappa_1 (\nZ_V(s) - \nZ_C(s))(\nK_2 -\nZ_C(s))\differential{s}\\
        &\quad  + \kappa_2^2\int_{0}^{t} n\left(F(\nZ_V(s), \nZ_C(s)- n^{-1}) - F(\nZ_V(s), \nZ_C(s)) \right)^2 \kappa_{-1}\nZ_C(s) \differential{s}\\
        =&\ \int_{0}^{t} \left(1 -n^{-1}\kappa_2 \lf(\SC{D}_{F,1}(\nZ_V(s), \nZ_C(s)-n^{-1})+\SC{D}_{F,2}(\nZ_V(s), \nZ_C(s))\ri)\right)^2 \\
        &\quad\quad  \times  \kappa_2 \nZ_C(s)\differential{s}\\
        &\quad + \kappa_2^2\int_{0}^{t} n^{-1}\lf(\SC{D}_{F,2}(\nZ_V(s), \nZ_C(s))\ri)^2  \kappa_1 (\nZ_V(s) - \nZ_C(s))  \\
        &\quad\quad \times (\nK_2 -\nZ_C(s))\differential{s}\\
         &\quad  + \kappa_2^2\int_{0}^{t}n^{-1}\lf(\SC{D}_{F,2}(\nZ_V(s), \nZ_C(s))\ri)^2 \kappa_{-1}\nZ_C(s) \differential{s}\\
         \equiv&\ \int_0^t \kappa_2 \nZ_C(s)\differential{s} + n^{-1} \SC{R}^{(n)}(t),
    \end{aligned}    
    \end{align}
    where the second equality is just a consequence of Taylor's expansion of first order and the terms $\SC{D}_{F,1}$ and $\SC{D}_{F,2}$, respectively, are first-order derivatives of $F$ with respect to first and second variables (evaluated at some intermediate points). Hence, using the estimates in \eqref{eq:poi-sol-2nd partial-1} and \eqref{eq:poi-sol-1st partial-1}, the process  $\SC{R}^{(n)}$ can be estimated as follows: for some constant $\tilde\const_{F}$,
    \begin{align*}
 \sup_{t\leq T}|\SC{R}^{(n)}(t)|\leq \tilde\const_{F}T(1+(\nK_1)^2) \eqstop 
    \end{align*}
 Since the sequence  $\{\nK_1 :  n\ge 1 \}$ is tight, it immediately follows that 
\begin{align}
     \label{eq:mart-rem-conv}
      n^{-1}\sup_{t\leq T}|\SC{R}^{(n)}(t)| \prt 0\eqcomma 
 \end{align}
 as $n\rt \infty$. The tightness of the process $\int_0^\cdot \nZ_C(s)\differential{s}$ in $C([0,T],\setOfPositiveReals)$ has already been established in \Cref{prop:rel_compactness}. It follows that the sequence of random variables $\predictableVariation{\mar^{(n)}}_T$ is tight in $\setOfPositiveReals$, and the sequence of processes $\predictableVariation{\mar^{(n)}}$ is tight in $C([0,T],\setOfPositiveReals)$. This immediately implies tightness of the sequence  $\{\sup_{t\leq T} |\mar^{(n)}(t)|\}$ in $\setOfPositiveReals$ by virtue of \Cref{lem:sup-mart-tight}, and the  $C$-tightness of $\mar^{(n)}$ in $D([0,T],\R)$ by \cite[Theorem 3.6]{Whitt2007MCLT}. 

 It is now clear from \eqref{eq:Un-SDE-0} that the tightness of the sequence of random variables $\{\sup_{t\leq T}|\nU_V(t)|\}$ is a consequence of \Cref{lem:tight-sup}, and  $C$-tightness of the sequence of processes $\{\nU_V\}$ in $D([0,T],\setOfPositiveReals)$ is a consequence of \Cref{cor:int-eq-tight} in Appendix~\ref{sec:aux-results}.
 
 Let $(U_V(0), \tilde K_2, U, Z_V, \occ)$ be a limit point of the sequence $(\nU_V(0), n^{1/2}(\nK_2-K_2), \nU_V, \nZ_V, \occ^{(n)})$ with the measure $\occ$ given by \eqref{eq:lim-occ-form-0}. Thus there exists a subsequence --- which by a slight abuse of notation, we continue to index by $n$ --- such that
  \begin{align}
     \label{eq:weak-conv-limit-clt}
      (\nU_V(0), n^{1/2}(\nK_2-K_2), \nU_V, \nZ_V, \occ^{(n)}) \nRT (U_V(0), \tilde K_2, U, Z_V, \occ).
 \end{align}
 


Identification of the stochastic process $U$ requires a second-order Taylor's expansion of the integral in the third term on the right side of \eqref{eq:U-exp-0}  giving
 \begin{align}
     \label{eq:eqb-expansion-2}
     \begin{aligned}
   & n^{1/2} \int_0^t  \lf(\eqb^{-}(\nZ_V(s)) - \eqb^{-}(Z_V(s))\ri) \differential{s} \\
   & =\ \int_0^t\partial\eqb^{-}(K_2,Z_V(s))\nU_V(s) \\
   & \quad  
   + n^{1/2}\int_0^t \SC{D}^{(2)}_{C,*}(\nZ_V(s),Z_V(s)) (\nZ_V(s)-Z_V(s))^2\differential{s}\\
   & = \int_0^t\partial\eqb^{-}(K_2,Z_V(s))\nU_V(s)\differential{s} \\
   &\quad 
   + n^{-1/2}\int_0^t \SC{D}^{(2)}_{C,*}(\nZ_V(s),Z_V(s)) (\nU_V(s))^2\differential{s}\eqcomma 
    \end{aligned}
 \end{align}
  where the expression of the first-order derivative of the mapping $z_V \mapsto \eqb^{-}(K_2,z_V)$ is given in \eqref{eq:deriv-eqb}, and $\SC{D}^{(2)}_{C,*}$ -- the term involving second order derivative of $\eqb^{-}(K_2,z_V)$ -- because of \eqref{eq:eqb-deriv-lin-2},  satisfies
\begin{align*}
    |\SC{D}^{(2)}_{C,*}(\nZ_V(s),Z_V(s))|\leq \const_{*,2}\lf(1+ 2\lf((\nK_1)^2+Z^2_V(s)\ri)\ri).
\end{align*}
Consequently, by the tightness of $\{\sup_{t\leq T}|\nU_V(s)|\}$, $\{\nK_1\}$, as $n\rt \infty$ it follows that 
\begin{align}\label{eq:deriv2-conv-0}
\begin{aligned}
  n^{-1/2}& \int_0^T |\SC{D}^{(2)}_{C,*}(\nZ_V(s),Z_V(s)) (\nU_V(s))^2|\differential{s}\\
  & \leq \ \const_{*,2}\lf(1+ 2\lf((\nK_1)^2+\sup_{t\leq T}Z^2_V(s)\ri)\ri)\sup_{t\leq T}|\nU_V(t)|^2T n^{-1/2}\\
 \ & \prt 0. 
\end{aligned}
\end{align}
 Plugging this in \eqref{eq:U-exp-0}, we get 
    \begin{align}\label{eq:Un-SDE-1}
    \begin{aligned}
        \nU_V(t) & = \nU_V(0) + n^{-1/2}\kappa_2\lf(F(\nZ_V(t), \nZ_C(t)) - F(\nZ_V(0), \nZ_C(0))\ri)\\
        &\quad +\kappa_2\int_0^t\partial\eqb^{-}(K_2,Z_V(s))\nU_V(s) \differential{s} + \mar^{(n)}(t)\\
        & \quad + n^{-1/2}\kappa_2\int_0^t \SC{D}^{(2)}_{C,*}(\nZ_V(s),Z_V(s)) (\nU_V(s))^2\differential{s} \\
        &\quad 
        -n^{1/2} \kappa_2\sum_{i=0}^2\int_{0}^{t} \err^{(n)}_{F,i}\lf(\nZ_V(s), \nZ_C(s)\ri)\differential{s}. 
    \end{aligned}
    \end{align}
Since $\occ^{(n)} \RT \occ$ as $n \to \infty$ (with $\occ$ given by \eqref{eq:lim-occ-form-0}), \eqref{eq:pred-quad-mart-F}, \eqref{eq:mart-rem-conv} and \Cref{lem:conv-occ-ZV} show that as $n\rt \infty$
  \begin{align*}
      \predictableVariation{\mar^{(n)}}_t \prt \kappa_2 \int_0^t\eqb^{-}(K_2,Z_V(s))\ \differential{s}\eqstop 
  \end{align*}
 where the convergence in probability holds as the limit is deterministic. 
    It is easy to see that for some constant $\tilde\const_{F,1}$
    \begin{align*}
        \Eof{\sup_{t\le T} \absolute{\mar^{(n)}(t) - \mar^{(n)}(t-)}} \le \tilde\const_{F,1} n^{-1} \EE(1+(\nK_1)^2) \nrt 0\eqstop
    \end{align*}
  
    Therefore, by the \ac{MCLT} \citep{Whitt2007MCLT,Ethier:1986:MPC}, we have 
    \begin{align}\label{eq:mclt}
        \mar^{(n)} \RT \mar \text{ where } \mar(t) \defeq \int_{0}^{t}\sqrt{{ \kappa_2}\eqb^{-}(K_2,Z_V(s))}\differential{W}(s)\eqstop
    \end{align}
    as $n\to \infty$, 
    where $W$ is a standard Wiener process. 

  Next, by the hypothesis on $\nK_2$ and \Cref{lem:conv-occ-ZV} and the continuous mapping theorem, it follows that
  \begin{align}\label{eq:E2-conv}
 \begin{aligned}
      & n^{1/2}\int_{0}^{t} \err^{(n)}_{F,2}\lf(\nZ_V(s), \nZ_C(s)\ri)\differential{s} \\
      &= \kappa_1\int_{\setOfPositiveReals\times \setOfPositiveReals \times [0, T]}(z_V-z_C)\partial_2F(z_V,z_C) \occ^{(n)}(\differential{z_C}\times \differential{z_V}\times \differential{s})\\
      &\quad  \times n^{1/2}(\nK_2-K_2)\\
      & = \int_{\setOfPositiveReals\times \setOfPositiveReals\times [0, T]}\f{z_V-z_C}{\kappa_1(\eqb^+(K_2,z_V)-z_C)} \occ^{(n)}(\differential{z_C}\times \differential{z_V}\times \differential{s})\\
      &\quad  \times n^{1/2}(\nK_2-K_2) \\
      &\quad \nrt \tilde K_2\int_0^t \f{Z_V(s)-\eqb^-(K_2,Z_V(s))}{\sqrt{D(K_2,Z_V(s))}} \differential{s},
  \end{aligned}
  \end{align}
  where $D(K_2,z_V)$ is given by \eqref{eq:disc}.
  
 Finally, we show that 
 \begin{align}
     \label{eq:weak-conv-lin}
    \int_0^\cdot \partial\eqb^{-}(K_2,Z_V(s))\nU_V(s)\differential{s} \stackrel{n\rt \infty}\RT \int_0^\cdot& \partial\eqb^{-}(K_2,Z_V(s))U_V(s))\differential{s}.
 \end{align}
 To this end, we first notice that because of the  $C$-tightness of the sequence $\{\nU_V : n \ge 1\}$, the limit point $U_V$ (see \eqref{eq:weak-conv-limit-clt}) almost surely has paths in $C([0,T],\R)$.
 Now by the Skorohod representation theorem, there exists a probability space $(\tilde \Omega, \tilde {F}, \tilde {P})$ and processes $\tilde{U}^{(n)}_V, \tilde {U}_V$ defined on this space such that 
\begin{align*}
   \tilde{U}^{(n)}_V \stackrel{d}= \nU_V, \quad  \tilde {U}_V \stackrel{d}= U_V, \quad \tilde{U}^{(n)}_V \nrt \tilde {U}_V, \quad \text{a.s. in } D([0,T], \R)
\end{align*}
Since $Z_V$ and $K_2$ are deterministic,
\begin{align}
    \label{eq:int-dist-equal}
    \begin{aligned}
   \int_0^\cdot \partial\eqb^{-}(K_2,Z_V(s))\tilde{U}^{(n)}_V(s)\differential{s}\ \stackrel{d} =&\ \int_0^\cdot \partial\eqb^{-}(K_2,Z_V(s))\nU_V(s)\differential{s}, \\
   \int_0^\cdot \partial\eqb^{-}(K_2,Z_V(s))\tilde U_V(s))\differential{s} \ \stackrel{d} =&\ \int_0^\cdot \partial\eqb^{-}(K_2,Z_V(s))U_V(s))\differential{s}
\end{aligned}
\end{align}
Clearly, $\tilde {U}_V \stackrel{d}= U_V$ implies that $\tilde {U}_V$ almost surely has paths in $C([0,T],\R)$.
 Therefore, by \cite[Chapter VI, Proposition 1.17]{jacod2003limit}, 
 \begin{align}
     \label{eq:Un-conv-U}
     \sup_{t\leq T}|\tilde{U}^{(n)}_V(t)-\tilde U_V(t)| \nrt 0, \quad \text{a.s.}
 \end{align} 
Consequently, 
 \begin{align*}
     \sup_{t\leq T} \Big|\int_0^t& \partial\eqb^{-}(K_2,Z_V(s))(\tilde{U}^{(n)}_V(s)- \tilde U_V(s))\differential{s}\Big|\\
     & \leq\   \sup_{t\leq T}|\partial\eqb^{-}(K_2,Z_V(t))|\sup_{t\leq T}|\tilde{U}^{(n)}_V(t)-\tilde U_V(t)|T\\
      & \leq\ \const_{*,1}(1+\sup_{t\leq T} Z_V(t))\sup_{t\leq T}|\tilde{U}^{(n)}_V(t)-\tilde U_V(t)|T \rt 0 \quad \text{a.s.},
 \end{align*}
 as $n \rt \infty$, which (because of \eqref{eq:int-dist-equal}) implies (actually equivalent to) \eqref{eq:weak-conv-lin}. 
Collectively, \eqref{eq:Un-SDE-1}, \eqref{eq:FE-conv-0}, \eqref{eq:deriv2-conv-0}, \eqref{eq:mclt}, \eqref{eq:E2-conv} and \eqref{eq:weak-conv-lin} establish that $U_V$ satisfies the \ac{SDE} in \eqref{eq:FCLT}, completing the proof of \Cref{thm:Z_V_FCLT}. 
    \end{proof}
    

\appendix

    \section{Auxiliary results}
    \label{sec:aux-results}
    In this section, we provide additional definitions, and certain auxiliary results used in the main text. 

    \begin{myDefinition}
        A collection of stochastic processes $\{\nU :n \ge 1\}$ is said to be $C$-tight in $D([0, T], \setOfReals)$ if the collection is relatively compact and hence tight in $D([0, T], \setOfReals)$, and limit points of every weakly convergent subsequence lie in the space $C([0, T], \setOfReals)$. 
    \end{myDefinition}

    For an element $x \in D([0, T], \setOfReals)$, define the modulus of continuity
    \begin{align}\label{eq:modu_cont_defn}
                    \modu(x,T,\delta) \defeq \sup_{\substack{t_1,t_2\in [0,T]\eqcomma \\ |t_1-t_1|\leq \delta}}|x(t_1)-x(t_2)|\eqstop 
    \end{align}

    \begin{myLemma}
        A collection of stochastic processes $\{\nU :n \ge 1\}$ with paths in $D([0, T], \setOfReals)$ is $C$-tight in $D([0, T], \setOfReals)$ if 
        the following two conditions are satisfied:
        \begin{enumerate}
            \item For each $t$ in a subset of $[0, T]$ that is dense in $[0, T]$, and that contains both $0$ and $T$, 
            \begin{align}
                \lim_{k\to \infty} \limsup_{n} \probOf{ \absolute{\nU(t)} \ge k} = 0\eqstop 
                \label{eq:c-tight-in-d-cond-1}
            \end{align}
            \item For $\vep>0, \ \eta>0$, there exist $0<\delta<1$ and $n_0\ge 0$ such that
            \begin{align}
               \sup_{n\geq n_0}\probOf{\modu(\nU,T,\delta)\geq \eta}\leq \vep\eqstop  
               \label{eq:c-tight-in-d-cond-2}
            \end{align}

        \end{enumerate}
            \label{lem:C-tight-in-D}
    \end{myLemma}
    \begin{proof}[Proof of \Cref{lem:C-tight-in-D}]
        The proof follows from \cite[Theorems 7.3 and 13.2]{Billingsley1999Convergence}, and the Corollary to Theorem 13.4 in \cite[pp. 142]{Billingsley1999Convergence}.
    \end{proof}

    \begin{myRemark}
        In the light of condition \eqref{eq:c-tight-in-d-cond-2} in \Cref{lem:C-tight-in-D}, the tightness of either sequence of real-valued random variables $\{\sup_{t \in [0, T]} \absolute{\nU(t)} : n\ge 1\}$ or the sequence $\{ \nU(0) : n\ge 1\}$ implies the condition \eqref{eq:c-tight-in-d-cond-1}. See \cite[Lemma 3.9]{Whitt2007MCLT}.
    \end{myRemark}
    \begin{myLemma}
        \label{lem:tight-sup}
     For each $n\geq 0$, let $\nU, A^{(n)}$ and $B^{(n)}$ be stochastic processes  satisfying
     $$\nU(t)=A^{(n)}(t)+\int_0^t B^{(n)}(s)\nU(s) \differential{s}.$$
     Assume that the sequences of $\setOfPositiveReals$-valued random variables $\{\sup_{t\leq T}|A^{(n)}(t)|\}$ and $\{\sup_{t\leq T}|B^{(n)}(t)|\}$ are tight. Then, $\{\sup_{t\leq T}|\nU(t)|\}$ is tight in $\setOfPositiveReals$.
    \end{myLemma}
    \begin{proof}[Proof of Lemma \ref{lem:tight-sup}]
        The proof readily follows from Gr\"onwall's inequality.
    \end{proof}

   \begin{myLemma}\label{lem:int-proc-tight}
      For each $n\geq 0$, let $\Phi^{(n)}$ and $\nY$ be stochastic processes  with paths in the  spaces $C([0,T],\R)$, $L^1([0,T],\R)$, respectively, satisfying
     $$\Phi^{(n)}(t)=\Phi^{(n)}(0)+\int_0^t \nY(s) \differential{s}.$$
     Assume that the sequences of $\setOfPositiveReals$-valued random variables $\{|\Phi^{(n)}(0)|\}$ and $\{\sup_{t\leq T}|\nY(t)|\}$ are tight. Then, the sequence $\{\Phi^{(n)}\}$ is tight in $C([0,T],\R).$
   \end{myLemma}

   \begin{proof}[Proof of Lemma \ref{lem:int-proc-tight}]
Since the sequence $\{|\Phi^{(n)}(0)|\}$ is tight, by \cite[Theorem 7.3]{Billingsley1999Convergence}, we need to show that for $\vep>0, \ \eta>0$, there exist $0<\delta<1$ and $n_0>0$ such that
\begin{align}
 \label{eq:modulus-of-cont}
   \sup_{n\geq n_0}\probOf{\modu(\Phi^{(n)},T,\delta)\geq \eta}\leq \vep
\end{align}
By tightness of $\{\sup_{t\leq T}|\nY(t)|\}$ in $\setOfPositiveReals$, find $R(\vep)$ such that for all $n>0$, $$\probOf{\sup_{t\leq T}|\nY(t)| > R(\vep)} \leq \vep.$$
Now clearly, $\modu(\Phi^{(n)},T,\delta) \leq \sup_{t\leq T}|\nY(t)|\delta.$
Choose $\delta>0$ small enough such that $\eta\delta^{-1} \geq R(\vep)$. It now follows that for all $n>0$,
\begin{align*}
    \probOf{\modu(\Phi^{(n)},T,\delta)\geq \eta} \leq&\ \probOf{\sup_{t\leq T}|\nY(t)| \geq \eta\delta^{-1}} \leq  \probOf{\sup_{t\leq T}|\nY(t)| \geq R(\vep)}\\
    \leq&\ \vep.
\end{align*}

   \end{proof}

\begin{myCorollary}\label{cor:int-eq-tight}
For each $n\geq 0$, let $\nU, A^{(n)}$ and $B^{(n)}$ be stochastic processes  satisfying
     $$\nU(t)=A^{(n)}(t)+\int_0^t B^{(n)}(s)\nU(s) \differential{s}.$$
Assume that the  $\{\sup_{t\leq T}|A^{(n)}(t)|\}$ and $\{\sup_{t\leq T}|B^{(n)}(t)|\}$ are tight in $\setOfPositiveReals$ and $\{A^{(n)}\}$ is tight in $D([0,T],\R)$. Then $\{\nU\}$ is tight in $D([0,T],\R)$. If $\{A^{(n)}\}$ is $C$-tight in $D([0,T],\R)$, then so is $\{\nU\}$. 
\end{myCorollary}

\begin{myLemma}
    \label{lem:sup-mart-tight}
 Let $\{\mar^{(n)}\}$ be a sequence of square integrable martingales such that $\{\<\mar^{(n)}\>_T\}$ is tight in $\setOfPositiveReals.$ Then, the sequence of random variables  $\sup_{t\leq T}|\mar^{(n)}(t)|$ is tight in $\setOfPositiveReals.$   
\end{myLemma}

\begin{proof}[Proof of Lemma \ref{lem:sup-mart-tight}]
    Let $\vep>0$. By the tightness of $\{\<\mar^{(n)}\>_T\}$, choose $R_1(\vep)$ such that 
    $$\sup_n \probOf{\<\mar^{(n)}\>_T>R_1(\vep)} \leq \vep/2.$$ Let $R_2(\vep) \equiv \vep^{-1/2} (2R_1(\vep))^{1/2}$. By Lenglart--Rebolledo inequality \cite[Lemma 3.7]{Whitt2007MCLT}, for all $n>0$,
    \begin{align*}
        \probOf{\sup_{t\leq T}|\mar^{(n)}(t)|>R_2(\vep)}\leq&\  R_1(\vep)/R_2^2(\vep)+ \probOf{\<\mar^{(n)}\>_T>R_1(\vep)}\\
        \leq &\ \vep/2+\vep/2=\vep.
    \end{align*}
\end{proof}

\section{Acronyms}
    
\begin{acronym}
	\acro{BDG}{Burkholder--Davis--Gundy}
	\acro{CDF}{Cumulative Distribution Function}
	\acro{CLT}{Central Limit Theorem}
	\acro{CM}{Configuration Model}
	\acro{CRN}{Chemical Reaction Network}
	\acro{CTMC}{Continuous Time Markov Chain}
	\acro{DSA}{Dynamic Survival Analysis}
	\acro{DTMC}{Discrete Time Markov Chain}
	\acro{ER}{Erd\"{o}s-R\'{e}nyi}
	\acro{ESI}{Enzyme-Substrate-Inhibitor}
	\acro{FCLT}{Functional Central Limit Theorem}
	\acro{FLLN}{Functional Law of Large Numbers}
	\acrodefplural{FLLN}[FLLNs]{Functional Laws of Large Numbers}
	\acro{FPT}{First Passage Time}
	\acro{GP}{Gaussian Process}
	\acrodefplural{GP}[GPs]{Gaussian Processes}
    \acro{HJB}{Hamilton Jacobi Bellman}
	\acro{iid}{independent and identically distributed}
	\acro{IPS}{Interacting Particle System}
	\acro{KL}{Kullback-Leibler}
	\acro{LDP}{Large Deviations Principle}
	\acro{LLN}{Law of Large Numbers}
	\acrodefplural{LLN}[LLNs]{Laws of Large Numbers}
	\acro{MCLT}{Martingale Central Limit Theorem}
	\acro{MCMC}{Markov Chain Monte Carlo}
	\acro{MFPT}{Mean First Passage Time}
	\acro{MGF}{Moment Generating Function}
	\acro{MLE}{Maximum Likelihood Estimate}
	\acro{MM}{Michaelis--Menten}
	\acro{ODE}{Ordinary Differential Equation}
	\acro{PDE}{Partial Differential Equation}
	\acro{PDF}{Probability Density Function}
	\acro{PGF}{Probability Generating Function}
	\acro{PMF}{Probability Mass Function}
	\acro{PRM}{Poisson Random Measure}
	\acro{psd}{positive semi-definite}
	\acro{QSSA}{Quasi-Steady State Approximation}
	\acro{rQSSA}{reversible QSSA}
    \acro{SDE}{Stochastic Differential Equation}
	\acro{SPDE}{Stochastic Partial Differential Equation}
	\acro{sQSSA}{standard QSSA}
	\acro{tQSSA}{total QSSA}
	\acro{whp}{with high probability}
\end{acronym}

\section*{Code and data availability}
We did not use any data for this study. 

\section*{Declaration of interest}
The authors declare no conflict of interest.

\section*{Declaration of generative AI in scientific writing}
During the preparation of this work the author(s) did not make use of any generative AI. 

\section*{Author contributions}
Both authors contributed equally to all aspects of work that led to this manuscript.

\section*{Funding}
Research of A. Ganguly is supported in part by  NSF DMS - 2246815 and Simons Foundation (via Travel Support for Mathematicians).

W. R. KhudaBukhsh was supported in part by an International Research Collaboration Fund (IRCF) awarded by the University of Nottingham, United Kingdom, a Scheme 4 `Research in Pairs' Grant 42360 by the London Mathematical Society (LMS)  and in part by the Engineering and Physical Sciences Research Council (EPSRC) [grant number EP/Y027795/1].

 \bibliographystyle{plain} 
 \bibliography{references}

\begin{thebibliography}{10}

\bibitem{Anderson:2011:CTM}
D.~F. Anderson and T.~G. Kurtz.
\newblock Continuous time markov chain models for chemical reaction networks.
\newblock In {\em Design and Analysis of Biomolecular Circuits}, pages 3--42.
  Springer, 2011.

\bibitem{Applebaum_2009Levy}
David Applebaum.
\newblock {\em L\'evy Processes and Stochastic Calculus}.
\newblock Cambridge Studies in Advanced Mathematics. Cambridge University
  Press, 2 edition, 2009.

\bibitem{ArKoNe06}
Vladimir~I. Arnold, Valery~V. Kozlov, and Anatoly~I. Neishtadt.
\newblock {\em Mathematical {A}spects of {C}lassical and {C}elestial
  {M}echanics}, volume~3 of {\em Encyclopedia of Mathematical Sciences}.
\newblock Springer-Verlag, Berlin, third edition, 2006.
\newblock [Dynamical systems. III], Translated from the Russian original by E.
  Khukhro.

\bibitem{Ball:2006:AAM}
K.~Ball, T.~G. Kurtz, L.~Popovic, and G.~A. Rempala.
\newblock Asymptotic analysis of multiscale approximations to reaction
  networks.
\newblock {\em Annals of Applied Probability}, 16(4):1925--1961, 2006.

\bibitem{Billingsley1999Convergence}
Patrick Billingsley.
\newblock {\em Convergence of Probability Measures}.
\newblock Wiley, 7 1999.

\bibitem{Bobrowski2005Functional}
A.~Bobrowski.
\newblock {\em Functional analysis for probability and stochastic processes}.
\newblock Cambridge University Press, Cambridge, 2005.
\newblock An introduction.

\bibitem{Borghans1996QSSA}
Jos{\'e}A.~M. Borghans, Rob~J. de~Boer, and Lee~A. Segel.
\newblock Extending the quasi-steady state approximation by changing variables.
\newblock {\em Bulletin of Mathematical Biology}, 58(1):43--63, 1996.

\bibitem{Bremaud2020PointProcess}
Pierre Br\'emaud.
\newblock {\em Point process calculus in time and space---an introduction with
  applications}, volume~98 of {\em Probability Theory and Stochastic
  Modelling}.
\newblock Springer, Cham, 2020.

\bibitem{Budhiraja2019Analysis}
Amarjit Budhiraja and Paul Dupuis.
\newblock {\em Analysis and Approximation of Rare Events: Representations and
  Weak Convergence Methods}.
\newblock Springer US, 2019.

\bibitem{BDG18}
Amarjit Budhiraja, Paul Dupuis, and Arnab Ganguly.
\newblock Large deviations for small noise diffusions in a fast markovian
  environment.
\newblock {\em Electron. J. Probab.}, 23:Paper No. 112, 33, 2018.

\bibitem{Choi2017BeyondMM}
Boseung Choi, Grzegorz~A. Rempala, and Jae~Kyoung Kim.
\newblock Beyond the michaelis--menten equation: Accurate and efficient
  estimation of enzyme kinetic parameters.
\newblock {\em Scientific Reports}, 7(1), 12 2017.

\bibitem{Crudu2012AAP}
A.~Crudu, A.~Debussche, A.~Muller, and O.~Radulescu.
\newblock Convergence of stochastic gene networks to hybrid piecewise
  deterministic processes.
\newblock {\em The Annals of Applied Probability}, 22(5), 10 2012.

\bibitem{Eilertsen2024Unreasonable}
Justin Eilertsen, Santiago Schnell, and Sebastian Walcher.
\newblock The unreasonable effectiveness of the total quasi-steady state
  approximation, and its limitations.
\newblock {\em J. Theoret. Biol.}, 583:Paper No. 111770, 8, 2024.

\bibitem{Ethier:1986:MPC}
S.~N. Ethier and T.~G. Kurtz.
\newblock {\em {Markov Processes: Characterization and Convergence}}, volume
  282.
\newblock John Wiley \& Wiley, 1986.

\bibitem{FeKaZa13}
Eugene~A. Feinberg, Pavlo~O. Kasyanov, and Nina~V. Zadoianchuk.
\newblock Fatou's lemma for weakly converging probabilities, 2013.

\bibitem{foll99}
Gerald~B Folland.
\newblock {\em Real analysis: modern techniques and their applications},
  volume~40.
\newblock John Wiley \& Sons, 1999.

\bibitem{FrWe08}
Mark Freidlin and Alexander Wentzell.
\newblock Some recent results on averaging principle.
\newblock In {\em Topics in stochastic analysis and nonparametric estimation},
  volume 145 of {\em IMA Vol. Math. Appl.}, pages 1--19. Springer, New York,
  2008.

\bibitem{FrWe12}
Mark~I. Freidlin and Alexander~D. Wentzell.
\newblock {\em Random perturbations of dynamical systems}, volume 260 of {\em
  Grundlehren der mathematischen Wissenschaften [Fundamental Principles of
  Mathematical Sciences]}.
\newblock Springer, Heidelberg, third edition, 2012.
\newblock Translated from the 1979 Russian original by Joseph Sz\"ucs.

\bibitem{GAK15}
Arnab Ganguly, Derya Altintan, and Heinz Koeppl.
\newblock Jump-diffusion approximation of stochastic reaction dynamics: error
  bounds and algorithms.
\newblock {\em Multiscale Model. Simul.}, 13(4):1390--1419, 2015.

\bibitem{BharatBaburModel}
Arnab Ganguly and Wasiur~R. KhudaBukhsh.
\newblock Functional limit theorems and parameter inference for multiscale
  stochastic models of enzyme kinetics, 2024.

\bibitem{Gaee2023Averaging}
Alexandre G\'enadot.
\newblock Averaging for slow-fast piecewise deterministic {M}arkov processes
  with an attractive boundary.
\newblock {\em J. Appl. Probab.}, 60(4):1439--1468, 2023.

\bibitem{HaRa02}
E.L. Haseltine and J.B. Rawlings.
\newblock Approximate simulation of coupled fast and slow reactions for
  stochastic chemical kinetics.
\newblock {\em Journal of Chemical Physics}, 117(15):6959--6969, 2002.

\bibitem{Khas68}
R.~Z. Hasminski{\u i}.
\newblock On the principle of averaging the {I}t\^o's stochastic differential
  equations.
\newblock {\em Kybernetika (Prague)}, 4:260--279, 1968.

\bibitem{HMS15}
Gang Huang, Michel Mandjes, and Peter Spreij.
\newblock Large deviations for {M}arkov-modulated diffusion processes with
  rapid switching.
\newblock {\em Stochastic Processes and their Applications}, 126(6):1785 --
  1818, 2016.

\bibitem{IkedaWatanabe2014Stochastic}
Nobuyuki Ikeda and Shinzo Watanabe, editors.
\newblock {\em {Stochastic Differential Equations and Diffusion Processes}}.
\newblock North Holland, 2014.

\bibitem{jacod2003limit}
Jean Jacod and Albert~N Shiryaev.
\newblock {\em Limit Theorems for Stochastic Processes}.
\newblock Springer-Verlag Berlin Heidelberg, 2003.

\bibitem{JaKr12}
T.~Jahnke and M.~Kreim.
\newblock Error bound for piecewise deterministic process modelling stochastic
  reaction systems.
\newblock {\em Multiscale Modeling and Simulation}, 10(4):1119--1147, 2012.

\bibitem{Kallenberg2017RandomMeasures}
Olav Kallenberg.
\newblock {\em Random measures, theory and applications}, volume~77 of {\em
  Probability Theory and Stochastic Modelling}.
\newblock Springer, Cham, 2017.

\bibitem{Kallenberg2021Foundations}
Olav Kallenberg.
\newblock {\em Foundations of modern probability}, volume~99 of {\em
  Probability Theory and Stochastic Modelling}.
\newblock Springer, Cham, third edition, 2021.

\bibitem{Kang:2013:STM}
H.-W. Kang and T.~G. Kurtz.
\newblock Separation of time-scales and model reduction for stochastic reaction
  networks.
\newblock {\em Annals of Applied Probability}, 23(2):529--583, 2013.

\bibitem{Kang:2014:CLT}
H.-W. Kang, T.~G. Kurtz, and L.~Popovic.
\newblock Central limit theorems and diffusion approximations for multiscale
  {M}arkov chain models.
\newblock {\em Annals of Applied Probability}, 24(2):721--759, 2014.

\bibitem{Kang2019QSSA}
Hye-Won Kang, Wasiur~R. KhudaBukhsh, Heinz Koeppl, and Grzegorz~A. Rempa{\l}a.
\newblock {Quasi-Steady-State Approximations Derived from the Stochastic Model
  of Enzyme Kinetics}.
\newblock {\em Bulletin of Mathematical Biology}, 81(5):1303--1336, 2019.

\bibitem{karatzas1991brownian}
Ioannis Karatzas and Steven~E Shreve.
\newblock {\em Brownian Motion and Stochastic Calculus}.
\newblock Springer-Verlag New York, 2 edition, 1991.

\bibitem{Kim2020Misuse}
Jae~Kyoung Kim and John~J. Tyson.
\newblock Misuse of the michaelis–menten rate law for protein interaction
  networks and its remedy.
\newblock {\em PLOS Computational Biology}, 16(10):e1008258, 10 2020.

\bibitem{Kim2014validity}
Jae Kyoung Kim, Krešimir Josić, and Matthew R. Bennett.
\newblock The validity of quasi-steady-state approximations in discrete
  stochastic simulations.
\newblock {\em Biophysical Journal}, 107(3):783–793, 8 2014.

\bibitem{Kurtz1992Averaging}
Thomas~G. Kurtz.
\newblock Averaging for martingale problems and stochastic approximation.
\newblock In Ioannis Karatzas and Daniel Ocone, editors, {\em Applied
  Stochastic Analysis}, pages 186--209, Berlin, Heidelberg, 1992. Springer
  Berlin Heidelberg.

\bibitem{Lip96}
Robert Liptser.
\newblock Large deviations for two scaled diffusions.
\newblock {\em Probab. Theory Related Fields}, 106(1):71--104, 1996.

\bibitem{Mao2024Averaging}
Yong-Hua Mao and Jinghai Shao.
\newblock Averaging principle for two time-scale regime-switching processes.
\newblock {\em Electron. J. Probab.}, 29:Paper No. 14, 21, 2024.

\bibitem{Neis90}
A.~Ne{\u\i}shtadt.
\newblock Averaging, capture into resonances, and chaos in nonlinear systems.
\newblock pages 261--273. Amer. Inst. Phys., New York, 1990.

\bibitem{Norris1997MarkovChains}
J.~R. Norris.
\newblock {\em Markov Chains}.
\newblock Cambridge University Press, 2 1997.

\bibitem{PaVe01}
E.~Pardoux and A.~Yu. Veretennikov.
\newblock On the {P}oisson equation and diffusion approximation. {I}.
\newblock {\em Ann. Probab.}, 29(3):1061--1085, 2001.

\bibitem{PaVe02}
\`E. Pardoux and A.~Yu. Veretennikov.
\newblock On {P}oisson equation and diffusion approximation. {II}.
\newblock {\em Ann. Probab.}, 31(3):1166--1192, 2003.

\bibitem{PaVe03}
\`E. Pardoux and A.~Yu. Veretennikov.
\newblock On {P}oisson equation and diffusion approximation. {II}.
\newblock {\em Ann. Probab.}, 31(3):1166--1192, 2003.

\bibitem{Pedersen2008tQSSA}
Morten~Gram Pedersen, Alberto~M. Bersani, Enrico Bersani, and Giuliana Cortese.
\newblock The total quasi-steady-state approximation for complex enzyme
  reactions.
\newblock {\em Mathematics and Computers in Simulation}, 79(4):1010–1019, 12
  2008.

\bibitem{RaAr03}
C.V. Rao and A.P. Arkin.
\newblock Stochastic chemical kinetics and the quasi-steady-state assumption:
  Application to the gillespie algorithm.
\newblock {\em Journal of Chemical Physics}, (11), 2003.

\bibitem{Schnell2000Enzymekinetics}
S~Schnell.
\newblock Enzyme kinetics at high enzyme concentration.
\newblock {\em Bulletin of Mathematical Biology}, 62(3):483–499, 5 2000.

\bibitem{Schnell2013Validity}
Santiago Schnell.
\newblock Validity of the michaelis–menten equation – steady‐state or
  reactant stationary assumption: that is the question.
\newblock {\em The FEBS Journal}, 281(2):464–472, 11 2013.

\bibitem{Segel:1975:EK}
I.~H. Segel.
\newblock {\em Enzyme kinetics}, volume 360.
\newblock Wiley, New York, 1975.

\bibitem{Segel:1988:VSS}
L.~A. Segel.
\newblock On the validity of the steady state assumption of enzyme kinetics.
\newblock {\em Bull. Math. Biol.}, 50(6):579--593, 1988.

\bibitem{Segel:1989:QSS}
L.~A. Segel and M.~Slemrod.
\newblock The quasi-steady-state assumption: a case study in perturbation.
\newblock {\em SIAM Rev.}, 31(3):446--477, 1989.

\bibitem{Shen2022Averaging}
Guangjun Shen, Wentao Xu, and Jiang-Lun Wu.
\newblock An averaging principle for stochastic differential delay equations
  driven by time-changed {L}\'evy noise.
\newblock {\em Acta Math. Sci. Ser. B (Engl. Ed.)}, 42(2):540--550, 2022.

\bibitem{Skorokhod1989Asymptotic}
A.~V. Skorokhod.
\newblock {\em Asymptotic methods in the theory of stochastic differential
  equations}, volume~78 of {\em Translations of Mathematical Monographs}.
\newblock American Mathematical Society, Providence, RI, 1989.
\newblock Translated from the Russian by H. H. McFaden.

\bibitem{Srivastava2025QSSA}
Kashvi Srivastava, Justin Eilertsen, Victoria Booth, and Santiago Schnell.
\newblock Accuracy versus predominance: Reassessing the validity of the
  quasi-steady-state approximation.
\newblock 2025.

\bibitem{Stiefenhofer:1998:QSS}
M.~Stiefenhofer.
\newblock Quasi-steady-state approximation for chemical reaction networks.
\newblock {\em J. Math. Biol.}, 36(6):593--609, 1998.

\bibitem{Tzafriri:2003:MMK}
A.~R. Tzafriri.
\newblock Michaelis-{M}enten kinetics at high enzyme concentrations.
\newblock {\em Bull. Math. Biol.}, 65(6):1111--1129, 2003.

\bibitem{Tzafriri:2004:TQS}
A.~R. Tzafriri and E.~R. Edelman.
\newblock The total quasi-steady-state approximation is valid for reversible
  enzyme kinetics.
\newblock {\em J. Theor. Biol.}, 226(3):303--313, 2004.

\bibitem{Tzafriri:2007:QSS}
A.~R. Tzafriri and E.~R. Edelman.
\newblock Quasi-steady-state kinetics at enzyme and substrate concentrations in
  excess of the {M}ichaelis--{M}enten constant.
\newblock {\em J. Theor. Biol.}, 245(4):737--748, 2007.

\bibitem{Veretennikov1990Averaging}
A.~Yu. Veretennikov.
\newblock On an averaging principle for systems of stochastic differential
  equations.
\newblock {\em Mat. Sb.}, 181(2):256--268, 1990.

\bibitem{Whitt2002StochLimits}
Ward Whitt.
\newblock {\em Stochastic-process limits}.
\newblock Springer Series in Operations Research. Springer-Verlag, New York,
  2002.
\newblock An introduction to stochastic-process limits and their application to
  queues.

\bibitem{Whitt2007MCLT}
Ward Whitt.
\newblock Proofs of the martingale {FCLT}.
\newblock {\em Probab. Surv.}, 4:268--302, 2007.

\bibitem{Wilkinson2018SMS}
Darren~J Wilkinson.
\newblock {\em Stochastic Modelling for Systems Biology}.
\newblock Chapman and Hall/CRC, 2018.

\bibitem{Xu2017Averaging}
Jie Xu.
\newblock {$L^p$}-strong convergence of the averaging principle for slow-fast
  {SPDE}s with jumps.
\newblock {\em J. Math. Anal. Appl.}, 445(1):342--373, 2017.

\end{thebibliography}








\end{document}